\documentclass[12pt]{article}
\usepackage{hyperref}
\usepackage{amsmath,amsthm,amssymb, amscd}
\usepackage{ytableau,}
\usepackage[all]{xy}
\usepackage{enumerate}

\renewcommand{\tilde}[1]{\widetilde{#1}}

\newtheorem{thm}{Theorem}[section]
\newtheorem{cor}[thm]{Corollary}
\newtheorem{lem}[thm]{Lemma}
\newtheorem{prop}[thm]{Proposition}
\newtheorem{dfn}[thm]{Definition}
\numberwithin{equation}{section}  

\begin{document}

\title{Crepant Resolutions of a Slodowy Slice in a Nilpotent Orbit Closure in $\mathfrak{sl}_N(\mathbf{C})$}
\footnotetext[0]{{\it 2010 Mathematics Subject Classification:} Primary 14J17, 14M15; Secondary 14E30, 14L24, 17B45.\\
{\it Keywords:} crepant resolutions, Slodowy slices, nilpotent orbits, quiver varieties.
}
\author{Ryo \textsc{Yamagishi}\footnote{Ryo Yamagishi: Departement of Mathematics, Graduate School of Science, Kyoto University,
Kita-shirakawa Oiwake-cho, Kyoto, 606-8502, JAPAN ;
email: ryo-yama@math.kyoto-u.ac.jp}}

\date{}

\maketitle

\begin{abstract}
One of our results of this article is that every (projective) crepant resolution of a Slodowy slice in a nilpotent orbit closure in $\mathfrak{sl}_N(\mathbf{C})$ can be obtained as the restriction of some crepant resolution of the nilpotent orbit closure. We also show that there is a decomposition of the Slodowy slice into other Slodowy slices with good properties. From this decomposition, one can count the number of crepant resolutions.
 \end{abstract}

\section{Introduction\label{intro}}

Throughout this article every algebraic variety is defined over $\mathbf{C}$, and when we speak of crepant resolutions, we only consider projective ones.

A normal (singular) variety $X$ is called a symplectic variety if it has a symplectic form $\omega$ on the regular part of $X$ and $\omega$ extends without poles to some resolution of $X$ (cf. [B]). It is an interesting object in birational geometry, representation theory and so on. A typical example of a symplectic variety is (the normalization of) a nilpotent orbit closure $\overline{O}$ in a complex simple Lie algebra $\mathfrak{g}$. It has a symplectic form $\omega_{KK}$ which is called the Kostant-Kirillov form on the regular part $O$. When we take an element $x$ from $\overline{O}$, we can consider an affine subspace $\mathcal{S}_x$ (defined in section \ref{sl} for $\mathfrak{g}=\mathfrak{sl}_{N}(\mathbf{C})$) of  $\mathfrak{g}$ which is called a Slodowy slice at $x$. Then the intersection $\overline{O}\cap\mathcal{S}_x$ is also a symplectic variety whose symplectic structure is given by the restriction of $\omega_{KK}$. When $\mathfrak{g}=\mathfrak{sl}_{N}(\mathbf{C})$, nilpotent orbits are parametrized by partitions of $N$. We let $O_d$ denote the nilpotent orbit in $\mathfrak{sl}_{N}(\mathbf{C})$ corresponding to a partition $d$ of $N$. If $x\in O_{d'}\subset\overline{O_d}$ for another partition $d'$ of $N$, the isomorphism class of $\mathcal{S}_x\cap\overline{O_d}$ depends only on $d$ and $d'$. So we let $\mathcal{S}_{d',d}$ denote the isomorphism class. In this article, we study $\mathcal{S}_{d',d}$ from the viewpoint of birational geometry. For the case where $\overline{O_d}$ is the nilpotent cone in $\mathfrak{sl}_{N}(\mathbf{C})$, $\mathcal{S}_{d',d}$ is studied in [LNS] by using a Poisson deformation.

To study the birational geometry of a symplectic variety $X$, we want to find a crepant resolution of $X$, which is a nonsingular relative minimal model of $X$. So it is natural to ask whether $X$ has a crepant resolution and how we can construct it if it exists. When $X=\overline{O_d}$, there always exist crepant resolutions and we can construct them explicitly. Let $a=[a_1,\dots,a_m]\in\mathbf{Z}_{>0}^m$ be the dual partition of $d$. For a permutation $\sigma\in\mathfrak{S}_m$, a flag of flag type $\sigma(a)=[a_{\sigma(1)},\dots,a_{\sigma(m)}]\in\mathbf{Z}_{>0}^m$ is a collection of vector subspaces $\{0\}=U_0\subset U_1\subset\cdots\subset U_{m-1}\subset U_m=\mathbf{C}^N$ such that $\mathrm{dim}\,U_i/U_{i-1}=a_{\sigma(i)}$ for $1\le i\le m$. Then $G=SL_N(\mathbf{C})$ acts transitively on $\mathcal{F}_{\sigma(a)}$, the set of flags of flag type $\sigma(a)$. Since the stabilizer group $P$ of a flag in $\mathcal{F}_{\sigma(a)}$ is a parabolic subgroup of $G$, $\mathcal{F}_{\sigma(a)}$ has the structure of a projective variety and we call it the flag variety of flag type $\sigma(a)$. The cotangent bundle $T^*\mathcal{F}_{\sigma(a)}$ is identified with the set $\{(x,U_\bullet)\in\overline{O_d}\times\mathcal{F}_a|x(U_{i+1})\subset U_i,\,0\le i\le m-1\}$. Then the first projection $$\{(x,U_
\bullet)\in\overline{O_d}\times\mathcal{F}_a|x(U_{i+1})\subset U_i,\,0\le i\le m-1\}\to\overline{O_d}$$
gives a crepant resolution. (This map can be identified with the Springer map $s_P:G\times^P\mathfrak{n}(P)\to\overline{O_d}\,;\,(g,x)\mapsto gxg^{-1}$; cf. Theorem \ref{thm:nam}).
By varying permutations $\sigma$, we get all crepant resolutions. Since different flag types give nonisomorphic crepant resolutions, the number of crepant resolutions of $\overline{O_d}$ is given by $N(d):=\sharp\{\sigma(a)\in\mathbf{Z}_{>0}^m|\sigma\in\mathfrak{S}_m\}$. 
We can also construct crepant resolutions of $\mathcal{S}_{d',d}$ from those of $\overline{O_d}$ by the following theorem.

\begin{thm}\label{thm:main1}
Let $\mathcal{S}_{d',d}$ be as above. Then the restriction of any crepant resolution of $\overline{O_d}$ to $\mathcal{S}_{d',d}$ is a crepant resolution of $\mathcal{S}_{d',d}$. Conversely, any crepant resolution of $\mathcal{S}_{d',d}$ can be obtained as the restriction of some crepant resolution of $\overline{O_d}$.
\end{thm}

The next natural question is counting how many crepant resolutions (up to isomorphisms) $\mathcal{S}_{d',d}$ has. This is answered by clarifying the chamber structure on the relative Picard group of a crepant resolution of $\mathcal{S}_{d',d}$. Now we explain this  briefly for a more general situation (see section \ref{resol} for definitions and details). Let $\pi:Y\to X$ be a crepant resolution of a rational Gorenstein singularity $X$. Then the $\mathbf{R}$-coefficient relative Picard group $\mathrm{Pic}(\pi)_{\mathbf{R}}$ is finite dimensional and it has the movable cone $\mathrm{Mov}(\pi)$ inside it. When we are given another birational projective morphism $\pi_i:Y_i\to X$ from a $\mathbf{Q}$-factorial variety $Y_i$ which is isomorphic to $Y$ in codimension one over $X$, the relative Picard group $\mathrm{Pic}(\pi_i)_{\mathbf{R}}$ is naturally isomorphic to $\mathrm{Pic}(\pi)_{\mathbf{R}}$. We can regard the $\pi_i$-ample cone in $\mathrm{Pic}(\pi_i)_{\mathbf{R}}$ as the cone in the movable cone $\mathrm{Mov}(\pi)$ via this isomorphism and we call it the ample cone of $\pi_i$ in  $\mathrm{Pic}(\pi)_{\mathbf{R}}$. It is known that there are finitely many ample cones, and the finiteness of (isomorphism classes of) crepant resolutions of $X$ follows from this fact. The closures of the ample cones cover $\mathrm{Mov}(\pi)$ and they form a chamber structure on $\mathrm{Mov}(\pi)$. Note that the natural isomorphism of relative Picard groups identifies movable cones and gives the same chamber structures on them.

If $X=\overline{O_d}$ or $\mathcal{S}_{d',d}$, then $\mathrm{Pic}(\pi)_{\mathbf{R}}$ is isomorphic to $H^2(Y,\mathbf{R})$, and all $\pi_i:Y_i\to X$ corresponding to the ample cones are crepant resolutions. Hence the set of crepant resolutions of $X$ has a natural one-to-one correspondence with the set of the ample cones. In general, this property holds if $X$ is an affine symplectic variety with a good $\mathbf{C}^*$-action [Nam3]. Moreover, by looking at the chamber structure, we can determine which pair of crepant resolutions are related by a flop. We have the following result.

\begin{thm}\label{thm:main2}
Let $Y$ be a crepant resolution of $\overline{O_d}$ and let $Y'$ be the restriction of $Y$ to $\mathcal{S}_{d',d}$. Then the restriction map $H^2(Y,\mathbf{R})\to H^2(Y',\mathbf{R})$ is surjective. Moreover, if it is an isomorphism, it maps the ample cones in $H^2(Y,\mathbf{R})$ to those in $H^2(Y',\mathbf{R})$ bijectively.
\end{thm}

The chamber structure of $H^2(Y,\mathbf{R})$ has already been studied in [Nam2] (cf. Theorem \ref{thm:nam}). By the above theorem, the chamber structure of $H^2(Y,\mathbf{R})$ coincides with that of $H^2(Y',\mathbf{R})$ if the restriction map is an isomorphism. In particular, the number of crepant resolutions of $\overline{O_d}$ is equal to that of $\mathcal{S}_{d',d}$ in this case. But in general, different resolutions of $\overline{O_d}$ may restrict to the same resolutions of $\mathcal{S}_{d',d}$. Even in this case, we can determine the chamber structure of $H^2(Y',\mathbf{R})$ by the following theorem.

\begin{thm}\label{thm:main3}
There is a decomposition $\mathcal{S}_{d',d}\cong\mathcal{S}_{d'^1,d^1}\times\cdots\times\mathcal{S}_{d'^r,d^r}$ ($r$ may be 1) into Slodowy slices in other $\mathfrak{sl}$'s such that the following properties hold.\\
(i) The product of crepant resolutions of $\mathcal{S}_{d'^i,d^i}$'s gives a crepant resolution of $\mathcal{S}_{d',d}$. Conversely every crepant resolution of $\mathcal{S}_{d',d}$ is of this form.\\
(ii) Let $Y_i$ be a crepant resolution of $\overline{O_{d^i}}$ and $Y'_i$ be the restriction of $Y_i$ to $\mathcal{S}_{d'^i,d^i}$. Then the restriction map $H^2(Y_i,\mathbf{R})\to H^2(Y'_i,\mathbf{R})$ is an isomorphism.\\
(iii) For $Y'=Y'_1\times\cdots\times Y'_r$, one has $H^2(Y',\mathbf{R})\cong\bigoplus_{i=1}^rH^2(Y'_i,\mathbf{R})$ and the set of ample cones in $H^2(Y',\mathbf{R})$ is equal to the set of products of ample cones in $H^2(Y_i',\mathbf{R})$ via this isomorphism.
\end{thm}

In section \ref{eg}, we illustrate how to get the decomposition above in terms of Young diagrams. The following is a direct corollary to the above theorems.

\begin{cor}
The number of crepant resolutions of $\mathcal{S}_{d',d}$ is $\prod_{i=1}^rN(d^i)$.
\end{cor}

We will prove the above theorems in section \ref{proof}. The main idea of the proofs is to realize $\mathcal{S}_{d',d}$ and its special crepant resolution as quiver varieties of type A via the isomorphisms constructed by Maffei [M2] (cf. Theorem \ref{thm:maf}). A quiver variety introduced by Nakajima [Nak1], [Nak2] is also a typical example of a symplectic variety (see section \ref{quiver}). It is obtained as a GIT quotient of some affine variety $\Lambda$ by some algebraic group $G_v$. To obtain a GIT quotient of $\Lambda$, one should take an element from the ($\mathbf{R}$-coefficient) character group $\chi(G_v)_{\mathbf{R}}$ of $G_v$. We say that two characters $\chi_1$ and $\chi_2$ are GIT equivalent if they give canonically isomorphic GIT quotients. By VGIT (Variation of Geometric Invariant Theory, cf. [DH] and [T]), the vector space $\chi(G_v)_{\mathbf{R}}$ has a chamber structure by GIT equivalence classes. We will construct an isomorphism between $\chi(G_v)_{\mathbf{R}}$ and $H^2(Y',\mathbf{R})$ for a crepant resolution  $Y'$ of  $\mathcal{S}_{d',d}$ so that some GIT chambers correspond to ample cones (Proposition \ref{prop:comm}). It also has the property that if a GIT chamber $C$ corresponds to an ample cone $A$, the GIT quotient associated to $C$ coincides with the crepant resolution of $\mathcal{S}_{d',d}$ associated to $A$. In particular, every crepant resolution of $\mathcal{S}_{d',d}$ can be obtained as a quiver variety.

Thanks to the above isomorphism between the cohomology group and the character group, we can interpret every crepant resolution as a GIT quotient (or a quiver variety). The decomposition of $\mathcal{S}_{d',d}$ in Theorem \ref{thm:main3} can be obtained easily from the quiver variety description. The decomposition  of $H^2(Y',\mathbf{R})$ (resp. the ample cones in $H^2(Y',\mathbf{R})$) in the theorem follows from the fact that the character group
$\chi(G_v)_{\mathbf{R}}$ (resp. the GIT chambers in $\chi(G_v)_{\mathbf{R}}$) decomposes into the character groups (resp. GIT chambers) associated to $Y_i'$'s.

\section{Nilpotent orbits and Slodowy slices in $\mathfrak{sl}_N(\mathbf{C})$\label{sl}}

In this section we recall the definitions and the basic properties of nilpotent orbits and Slodowy slices in $\mathfrak{sl}_N(\mathbf{C})$. Consider the special linear group $SL_N(\mathbf{C})$ and its Lie algebra $\mathfrak{sl}_N(\mathbf{C})=\{N\times N\mathrm{-matrix}\;A|\mathrm{tr}\,A=0\}$. $SL_N(\mathbf{C})$ acts on $\mathfrak{sl}_N(\mathbf{C})$ by conjugation and we call the orbits of this action {\em adjoint orbits}.  Since conjugation preservers nilpotency, every adjoint orbit of a nilpotent matrix in $\mathfrak{sl}_N(\mathbf{C})$ consists of nilpotent matrices. We call such orbits {\em nilpotent adjoint orbits} or simply {\em nilpotent orbits}. Since $SL_N(\mathbf{C})$-conjugacy classes coincide with $GL_N(\mathbf{C})$-conjugacy classes, nilpotent orbits in $\mathfrak{sl}_N(\mathbf{C})$ are in one-to-one correspondence with Jordan normal forms whose diagonals are zero. So nilpotent orbits also have one-to-one correspondence with partitions of $N$, where a {\em partition} of $N$ is a tuple $d=[d_1,\dots,d_l]$ of natural numbers $d_i$ such that $d_1\ge\cdots\ge d_l$ and $\sum_{i=1}^ld_i=N$. In this article we often identify a partition $d$ with a Young diagram whose $i$-th row has $d_i$ boxes. For example, $d=[4,4,2,1]$ is identified with the following Young diagram.
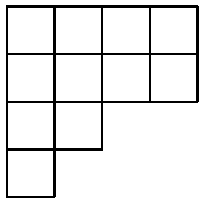
\begin{center}
\ydiagram{4,4,2,1}
\end{center}

For a partition $d$ of $N$, let $O_d$ denote the corresponding nilpotent orbit. The set of nilpotent orbits has a partial order $\le$ defined as follows:
$$O_{d'}\le O_d\stackrel{\mathrm{def}}{\iff}\overline{O_{d'}}\subset\overline{O_d}$$
where the overlines mean closures in $\mathfrak{sl}_N(\mathbf{C})$.
The following proposition due to Gerstenhaber [CM, 5.1] enables us to decide when $O_{d'}\le O_d$ occurs in terms of the partitions $d,d'$.

\begin{prop}
Let $d=[d_1,\cdots,d_l]$ and $d'=[d'_1,\cdots,d'_m]$ be two partitions of $N$. Then $O_{d'}\le O_d$ if and only if
$$\sum_{i=1}^kd'_i\le\sum_{i=1}^kd_i\hspace{5mm}for\;all\;1\le k\le\mathrm{min}\{l,m\}.$$
\end{prop}

Next we define Slodowy slices, which were introduced in [S]. Let $x$ be a nilpotent element in $\mathfrak{sl}_N(\mathbf{C})$. By the Jacobson-Morozov Theorem, there are two elements $y$ and $h$ in $\mathfrak{sl}_2(\mathbf{C})$ such that
$$[x,y]=h, \;\;[h,x]=2x \;\text{ and } \;[h,y]=-2y$$
(such a triple ($x,y,h$) is called an $\mathfrak{sl}_2(\mathbf{C})$-triple if $x\ne0$). A {\em Slodowy slice} at $x$ is the affine subspace $$\mathcal{S}_x:=x+\mathrm{Ker\,(ad}\,y)$$
of $\mathfrak{sl}_N(\mathbf{C})$ where $\mathrm{ad}\,y$ is the $\mathbf{C}$-linear map $\mathfrak{sl}_N(\mathbf{C})\to\mathfrak{sl}_N(\mathbf{C})$ defined by $v\mapsto[y,v]$. Note that $\mathfrak{sl}_2(\mathbf{C})$-triples are not unique but they are all $SL_N(\mathbf{C})$-conjugate. Let $O_{d'}$ be the nilpotent orbit of $x$ where $d'$ is a partition of $N$. Slodowy slices have the following properties ([CG, 3.7]):\\
$\bullet\;\mathcal{S}_x\cap O_{d'}=\{x\}$.\\
$\bullet\;\mathcal{S}_x$ intersects $O_{d'}$ transeversally, i.e.
$$T_x\mathcal{S}_x\oplus T_xO_{d'}=\mathfrak{sl}_N(\mathbf{C}).$$

Let $d$ and $d'$ be two partitions of $N$, and $x$ be a nilpotent element in $O_{d'}$. We assume that $x$ is in $\overline{O_d}$ or equivalently $O_{d'}\le O_d$. 
Since $\mathcal{S}_x\cap\overline{O_d}$ and $\mathcal{S}_{x'}\cap\overline{O_d}$ are isomorphic if $x'\in O_{d'}$, we define $\mathcal{S}_{d',d}$ as the isomorphism class of $\mathcal{S}_x\cap\overline{O_d}$. Note that if $x=0$ then $\mathcal{S}_{d',d}=\overline{O_d}$. Therefore we can think of $\mathcal{S}_{d',d}$'s as a generalization of nilpotent orbit closures. One important feature of $\mathcal{S}_{d',d}$ is that $\mathcal{S}_{d',d}$ is a symplectic variety (or has a symplectic singularity) (cf. [B]).

\begin{dfn}
A normal algebraic variety $X$ is a $\textbf{symplectic\;variety}$ if the regular locus of $X$ admits symplectic form $\omega$ and for some resolution $\pi:Y\to X$, the pullback $\pi^*(\omega)$ extends to whole $Y$ without poles.
\end{dfn}\vspace{3mm}

It is known that $\overline{O_d}$ is normal [KP] and its regular locus ${O_d}$ has a symplectic form called the {\em Kostant-Kirillov form} [CG, 1.1]. As we will see in section \ref{quiver}, the singular variety $\overline{O_d}$ always has a resolution from a symplectic manifold, and hence the following proposition holds. 

\begin{prop}
A nilpotent orbit closure $\overline{O_d}$ in $\mathfrak{sl}_N(\mathbf{C})$ is a symplectic variety.
\end{prop}

The restriction of the Kostant-Kirillov form on $O_d$ to $\mathcal{S}_x\cap O_d$ gives a symplectic form on $\mathcal{S}_x\cap O_d$. We know from the argument in [S, Chs.4,5] that a resolution of $\mathcal{S}_{d',d}$ can be obtained by restricting some resolution of $\overline{O_d}$. From these facts we conclude that $\mathcal{S}_{d',d}$ is a symplectic variety as well.\vspace{3mm}

\noindent{\bf Remark.} 1. We can define nilpotent orbits and Slodowy slices for other complex simple Lie algebras in the same way. But nilpotent orbit closures in other complex simple Lie algebras are not normal in general. In this case, their normalizations are symplectic varieties [P].

2. For $\mathfrak{sl}_N(\mathbf{C})$, the fact that $\mathcal{S}_{d',d}$ is a symplectic variety also follows from a quiver variety description (see section \ref{quiver}).

\section{Birational geometry of crepant resolutions\label{resol}}
Our main theorems are concerned with birational geometry of crepant resolutions of $\mathcal{S}_{d',d}$. In this section we introduce several notions for resolutions of more general singularities.

\begin{dfn}
Let $X$ be a normal (singular) variety.\\
$\bullet\;X$ has a $\textbf{rational singularity}$ if $R^i\pi_*\mathcal{O}_Y=0$ for every resolution $\pi:Y\to X$ and all \text{ \;}$i>0$.\\
$\bullet\;X$ is $\textbf{Gorenstein}$ if the canonical divisor $K_X$ is Cartier. \\
$\bullet$ A resolution $\pi:Y\to X$ of $X$ is $\textbf{crepant}$ if $K_Y=\pi^*K_X$.
\end{dfn}

Let  $\pi:Y\to X$ be a crepant resolution of a rational Gorenstein singularity $X$. We can define the {\em relative Picard group} of $\pi$ as $\mathrm{Pic}(\pi):=\mathrm{Pic}(Y)/\mathrm{Im}\,\pi^*$ where $\pi^*:\mathrm{Pic}(X)\to\mathrm{Pic}(Y)$ is the pullback homomorphism. From now on, for any Abelian group $A$, let $A_\mathbf{R}$ denote the $\mathbf{R}$-vector space $A\otimes_\mathbf{Z}\mathbf{R}$. The rationality of the singularity of $X$ implies the following.

\begin{lem}\label{lem:pic}
$\mathrm{Pic}(\pi)_\mathbf{R}$ is a finite dimensional $\mathbf{R}$-vector space. If $X$ is affine and contractible, then $\mathrm{Pic}(\pi),\mathrm{Pic}(Y)$ and $H^2(Y,\mathbf{Z})$ are isomorphic to each other.
\end{lem}
{\em Proof.} By the Leray spectral sequence, there is the following exact sequence of Abelian groups:
\begin{equation}
0\to H^1(X,\mathcal{O}_X^*)\xrightarrow{\pi^*} H^1(Y,\mathcal{O}_Y^*)\to H^0(X,R^1\pi_*\mathcal{O}_Y^*)\to H^2(X,\mathcal{O}_X^*).
\end{equation}
So $\mathrm{Pic}(\pi)=H^1(Y,\mathcal{O}_Y^*)/\pi^*H^1(X,\mathcal{O}_X^*)$ is a subgroup of $H^0(X,R^1\pi_*\mathcal{O}_Y^*)$. The exponential exact sequence
$$0\to\mathbf{Z}\to\mathcal{O}_Y\to \mathcal{O}_Y^*\to1$$
induces a long exact sequence
\begin{equation}
\begin{aligned}0&\to\pi_*\mathbf{Z}\to\pi_*\mathcal{O}_Y\to\pi_*\mathcal{O}_Y^*\\
&\to R^1\pi_*\mathbf{Z}\to R^1\pi_*\mathcal{O}_Y\to R^1\pi_*\mathcal{O}_Y^*\\
&\to R^2\pi_*\mathbf{Z}\to R^2\pi^*\mathcal{O}_Y\to\cdots.\end{aligned}
\end{equation}
As $R^1\pi_*\mathcal{O}_Y=R^2\pi_*\mathcal{O}_Y=0$, the induced map $H^0(X,R^1\pi_*\mathcal{O}_Y^*)\to H^0(X,R^2\pi_*\mathbf{Z})$ is an isomorphism. Finite generation of the relative N\'{e}ron-Severi group $H^0(X,R^2\pi_*\mathbf{Z})$ implies that of $\mathrm{Pic}(\pi)$. 

Suppose that $X$ is affine and contractible. Then we have $\mathrm{Pic}(Y)=\mathrm{Pic}(\pi)=H^0(X,R^2\pi_*\mathbf{Z})$ since $H^1(X,\mathcal{O}_X^*)=H^2(X,\mathcal{O}_X^*)=0$ in (3.1). On the other hand, by the Leray spectral sequence again, there is a filtration $H^2(Y,\mathbf{Z})=F^0\supset F^1\supset F^2$ by Abelian groups such that
$$\begin{aligned}F^0/F^1&=\mathrm{Ker}\,(H^0(X,R^2\pi_*\mathbf{Z})\to H^2(X,R^1\pi_*\mathbf{Z}))\\
F^1/F^2&=\mathrm{Ker}\,(H^1(X,R^1\pi_*\mathbf{Z})\to H^3(X,\mathbf{Z}))\\
F^2&=\mathrm{Coker}\,(H^0(X,R^1\pi_*\mathbf{Z})\to H^2(X,\mathbf{Z})).\end{aligned}$$
Since $\pi_*\mathcal{O}_Y\to\pi_*\mathcal{O}_Y^*$ in (3.2) coincides with the exponential map $\mathcal{O}_X\to\mathcal{O}_X^*$, and $R^1\pi_*\mathcal{O}_Y=0$, we have $R^1\pi_*\mathbf{Z}=0$. Therefore $F^1=F^2=0$ and $F^0=H^0(X,R^2\pi_*\mathbf{Z})$. \qed

\vspace{3mm}
The following lemmas are important when one treats all crepant resolutions.

\begin{lem}\label{lem:cre}
If $\pi':Y'\to X$ is another crepant resolution of $X$, then $Y$ and $Y'$ are isomorphic in codimension one over $X$.
\end{lem}
{\em Proof.} Let $W'$ be the irreducible component of $Y\times_X Y'$ which dominates $Y$ and $Y'$, and let $W$ be a resolution of $W'$. Then there is the following commutative diagram where $r_1:W\to Y$ and $r_2:W\to Y'$ are both birational morphisms of smooth varieties.

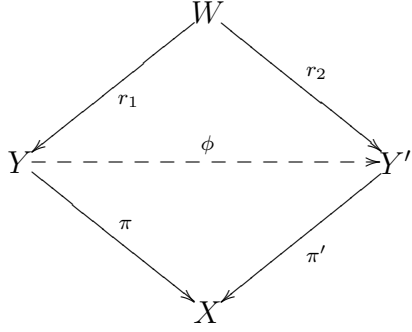
\vspace{0.5cm}
\begin{xy}
(0,0)*{}="",(50,20)*{Y}="A",(100,20)*{Y'}="B",(75,0)*{X}="C",(75,40)*{W}="D"
\ar @{-->}^{\phi}"A";"B"
\ar^{\pi} "A";"C"
\ar^{\pi'} "B";"C"
\ar^{r_1} "D";"A"
\ar^{r_2} "D";"B"
\end{xy}
\vspace{0.5cm}

We can write $K_W=r_1^*(K_Y)+E_1=r_2^*(K_Y)+E_2$ where $E_i$ is a full $r_i$-exceptional divisor. Since $\pi$ and $\pi'$ are both crepant, we have $E_1=E_2$. We must show that $\phi(D)$ is a divisor of $Y'$ for every divisor $D$ of $Y$. Otherwise $r_1^*(D)$ would be $r_2$-exceptional but not $r_1$-exceptional. This is a contradiction. \qed

\begin{lem}\label{lem:iso}
Let $\pi':Y'\to X$ be a projective morphism from a $\mathbf{Q}$-factorial variety $Y'$ and let $L'$ be a $\pi'$-ample line bundle on $Y'$. If the rational map $f=\pi'^{-1}\circ\pi:Y\dashrightarrow Y'$ is an isomorphism in codimension one, then there is a natural isomorphism of $\mathbf{R}$-vector spaces $f^*:\mathrm{Pic}(\pi')_\mathbf{R}\cong\mathrm{Pic}(\pi)_\mathbf{R}$, and $Y'$ is isomorphic to $\mathrm{\mathbf{Proj}}_X\bigoplus_{i=0}^\infty \pi_*(f^*L'^{\otimes i})$.
\end{lem}
{\em Proof.} Since $Y$ and $Y'$ are $\mathbf{Q}$-factorial, the first claim is clear. The second claim follows since $\pi_*(f^*L'^{\otimes i})\cong\pi'_*(L'^{\otimes i})$ for all $i$.\qed

\begin{dfn}
A line bundle $L$ on $Y$ is $\boldsymbol{\pi}$\textbf{-movable} if $\mathrm{codim\,Supp(Coker}\,\alpha)\ge2$ where $\alpha:\pi^*\pi_*L\to L$ is the natural map of sheaves on $Y$. The $\boldsymbol{\pi}$\textbf{-movable cone} $\mathrm{Mov}(\pi)$ in $\mathrm{Pic}(\pi)_\mathbf{R}$ is the cone generated by the classes of $\pi$-movable line bundles.
\end{dfn}

Now we assume that $X$ is affine. By the result of [BCHM], $\pi$ is a relative version of a Mori dream space (cf. [HK]). So $\pi$ has following properties:\\
There are finitely many projective morphisms $\{\pi_i:Y_i\to X\}_i$ from $\mathbf{Q}$-factorial varieties with $f_i:=\pi_i^{-1}\circ\pi$ which are isomorphisms in codimension one over $X$ such that:\\
(1) $f_i^*\mathrm{Amp}(\pi_i)$'s are disjoint,\\
(2) $\mathrm{Mov(\pi)}=\bigcup_if_i^*\overline{\mathrm{Amp}}(\pi_i)$,\\
(3) $\mathrm{Mov(\pi)}$ and $f_i^*\overline{\mathrm{Amp}}(\pi_i)$'s are polyhedral cones in $\mathrm{Pic}(\pi)_\mathbf{R}$, and\\
(4) $\pi_i$ and $\pi_j$ are related by a flop if and only if the cones $f_i^*\overline{\mathrm{Amp}}(\pi_i)$ and $f_j^*\overline{\mathrm{Amp}}(\pi_j)$ are adjacent in $\mathrm{Pic}(\pi)_\mathbf{R}$.\\
We say that $f_i^*\mathrm{Amp}(\pi_i)$ is the {\em ample cone} of $\pi_i$ in $\mathrm{Pic}(\pi)_\mathbf{R}$. Note that $Y_i$ and $Y_j$ are not isomorphic over $X$ if $i\ne j$. Indeed, otherwise $f_j^{-1}\circ f_i:Y_i\dashrightarrow Y_j$ would extend to be an isomorphism over $X$, contrary to the fact that $f_i^*\mathrm{Amp}(\pi_i)$'s are disjoint.

Note that a nilpotent orbit closure $\overline{O_d}$ in $\mathfrak{sl}_N(\mathbf{C})$ has a rational Gorenstein singularity [B; 1.3]. To end this section, we describe all crepant resolutions of $\overline{O_d}$. The results about crepant resolutions of $\overline{O_d}$ (which are not necessary to prove the theorems in section \ref{intro}) are summarized as follows.

\begin{thm}\label{thm:nam}([Nam1] and [Nam2] for $\mathfrak{g}=\mathfrak{sl}_N(\mathbf{C})$)
Let $G=SL_N(\mathbf{C})$. Then:\\
(a) There is a parabolic subgroup $P$ of $G$ such that the Springer map $s_P:G\times^P\mathfrak{n}(P)\to\overline{O_d}$ defined by $(g,x)\mapsto gxg^{-1}$ gives a crepant resolution of $\overline{O_d}$ where $\mathfrak{n}(P)$ denotes the nilradical of $\mathrm{Lie}(P)$.\\
(b) For a parabolic subgroup $P_i$ of $G$ which has the same Levi part $L$ of $P$, the Springer map $s_{P_i}:G\times^{P_i}\mathfrak{n}(P_i)\to \overline{O_d}$ also gives a crepant resolution of $\overline{O_d}$.\\
(c) Every crepant resolution of $\overline{O_d}$ is of the form in (b) up to isomorphisms. \\
(d) For two parabolic subgroups $P_i$ and $P_j$ of $G$, the Springer maps $s_{P_i}$ and $s_{P_j}$ give the isomorphic resolutions of $\overline{O_d}$ if and only if $P_i$ and $P_j$ are $G$-conjugate.\\
(e) The symmetric group $\mathfrak{S}_m$ for some $m$ naturally acts on $\mathrm{Pic}(G\times^P\mathfrak{n}(P))_{\mathbf{R}}$ and the closure of the ample cone is the fundamental domain of this action.\\
(f) There is a finite group $W_L$ associated to $L$ such that $W_L$ naturally acts on $\mathrm{Pic}(G\times^P$ $\mathfrak{n}(P))_{\mathbf{R}}$ and the movable cone is the fundamental domain of this action.
\end{thm}

In the next section we will describe the crepant resolutions $G\times^{P_i}\mathfrak{n}(P_i)$ as the cotangent bundles of flag varieties.

\section{Quiver varieties of type A and GIT chambers\label{quiver}}
In this section we introduce quiver varieties of type A and prepare some GIT notions for the next section. We will define quiver varieties not as hyper-K\"{a}hler quotients but as GIT quotients, as in [Nak2].
Let $n$ be a natural number. Take complex vector spaces $V_i$ and $W_i$ for each $i=1,\cdots,n$ whose dimensions are $v_i$ and $w_i$, which may be zero. Set $v=(v_1,\cdots,v_n)$ and $w=(w_1,\cdots,w_n)$. We define a {\em double quiver representation}  associated with the dimension vectors $v$ and $w$ as
$$M(v,w)=\bigoplus_{i=1}^{n-1}(\mathrm{Hom}(V_i,V_{i+1})\oplus\mathrm{Hom}(V_{i+1},V_i))\oplus\bigoplus_{j=1}^n(\mathrm{Hom}(W_j,V_j)\oplus\mathrm{Hom}(V_j,W_j))$$ where $\mathrm{Hom}$ means the vector space of all $\mathbf{C}$-linear maps. An element of $M(v,w)$ is usually written as a collection $(A_i,B_i, \Gamma_j,\Delta_j)_{1\le i\le n-1,1\le j\le n}$ where $A_i,\,B_i,\,\Gamma_j$ and $\Delta_j$ are elements of $\mathrm{Hom}(V_i,V_{i+1}),\,\mathrm{Hom}(V_{i+1},V_i),\,\mathrm{Hom}(W_j,V_j)$ and $\mathrm{Hom}(V_j,W_j)$ respectively. An element $(A_i,B_i, \Gamma_j,\Delta_j)_{i,j}$ of $M(v,w)$ can be understood as the configuration of linear maps in the diagram below.
\vspace{3mm}

\begin{xy}
(0,0)*{}="",(40,0)*{V_1}="1",(60,0)*{V_2}="2",(80,0)*{\;\cdots\;}="dot",(100,0)*{V_{n-1}}="n-1",(120,0)*{V_n}="n",(40,20)*{W_1}="w1",(60,20)*{W_2}="w2",(100,20)*{W_{n-1}}="wn-1",(120,20)*{W_n}="wn"
\ar @<1mm>"1";"2" ^{A_1} \ar @<1mm> "2";"dot"^{A_2} \ar @<1mm> "dot";"n-1"^{A_{n-2}} \ar @<1mm> "n-1";"n"^{A_{n-1}}
\ar @<1mm> "2";"1"^{B_1} \ar @<1mm> "dot";"2"^{B_2} \ar @<1mm> "n-1";"dot"^{B_{n-2}} \ar @<1mm> "n";"n-1"^{B_{n-1}}
\ar @<1mm>"w1";"1"^{\Gamma_1} \ar @<1mm>"w2";"2"^{\Gamma_2} \ar @<1mm>"wn-1";"n-1"^{\Gamma_{n-1}} \ar @<1mm>"wn";"n"^{\Gamma_n} \ar @<1mm>"1";"w1"^{\Delta_1} \ar @<1mm>"2";"w2"^{\Delta_2} \ar @<1mm>"n-1";"wn-1"^{\Delta_{n-1}} \ar @<1mm>"n";"wn"^{\Delta_n}
\end{xy}
\vspace{3mm}

Let $G_v=\prod_{i=1}^nGL(V_i)$ be the product of the general linear groups. It acts on $M(v,w)$ as
$(g_1,\cdots,g_n):(A_i,B_i, \Gamma_j,\Delta_j)_{i,j}\mapsto(g_{i+1}A_ig_i^{-1},g_iB_ig_{i+1}^{-1}, g_j\Gamma_j,\Delta_jg_j^{-1})_{i,j}$ for $g_i\in GL(V_i)$. The vector space $M(v,w)$ has a natural symplectic structure which is preserved by the action of $G_v$. The associated moment map is given by
$$\mu:M(v,w)\to\bigoplus_{i=1}^n\mathfrak{gl}(V_i),$$
$$\begin{aligned}\mu(A_i,B_i, \Gamma_j,\Delta_j)_{i,j}=&(\Gamma_1\Delta_1-B_1A_1\\&,\Gamma_2\Delta_2+A_1B_1-B_2A_2\\&\hspace{20mm}\vdots\\
&,\Gamma_{n-1}\Delta_{n-1}+A_{n-2}B_{n-2}-B_{n-1}A_{n-1}\\&,\Gamma_n\Delta_n+A_{n-1}B_{n-1}).\end{aligned}$$
Note that $\mu$ is $G_v$-equivariant where $G_v$ acts on $\bigoplus_{i=1}^n\mathfrak{gl}(V_i)$ by componentwise conjugation. Hence $\Lambda=\Lambda(v,w):=\mu^{-1}(0)$ is $G_v$-invariant. Let $\chi(G_v):=\mathrm{Hom}(G_v,\mathbf{C}^*)$ be the group  of characters of $G_v$. Since $\mathbf{Z}\to\mathrm{Hom}(GL(V_i),\mathbf{C}^*), m\mapsto\mathrm{det}^m$ is an isomorphism for nonzero $V_i$, $\chi(G_v)$ is isomorphic to $\mathbf{Z}^{n_0}$ where $n_0$ is the number of nonzero $V_i$'s. We sometimes identify an $n_0$-tuple of integers with a character of $G_v$ via this isomorphism.
For a character $\chi$ of $G_v$, set
$$R(\chi)=\{f\in\Gamma(\Lambda,\mathcal{O}_\Lambda)|g\cdot f=\chi(g)f\mbox{ for all }g\in G_v\}.$$
We say that a point $x\in\Lambda$ is $\chi$-{\em semistable} if there exist $i\in\mathbf{Z}_{>0}$ and $f\in R(\chi^i)$ such that $f(x)\ne0$. If moreover $x$ has a finite stabilizer and the $G_v$-orbit of $x$ is closed in $\{x\in\Lambda|f(x)\ne0\}$, we say that $x$ is $\chi$-{\em stable}. Let $\Lambda_\chi^{ss}$ (resp. $\Lambda_\chi^s$) denote the subset of $\chi$-semistable (resp. $\chi$-stable) points in $\Lambda$. We define the {\em quiver variety} of type A associated with $v,w$ and $\chi$ as
$$\mathfrak{M}_\chi(v,w):=\Lambda/\!/_\chi G_v:=\mathrm{Proj}\bigoplus_{i=0}^\infty R(\chi^i).$$

By definition, there is a projective morphism $\pi_\chi:\mathfrak{M}_\chi(v,w)\to\mathfrak{M}_\mathbf{0}(v,w)=\mathrm{Spec}\,R(\mathbf{0})$ where $\mathbf{0}=(0,\cdots,0)\in\mathbf{Z}^{n_0}\cong\chi(G_v)$ be the trivial character. The following is a fundamental result in GIT (cf. [MFK]).

\begin{prop}\label{prop:git}
The morphism $q:\Lambda_\chi^{ss}\to\Lambda/\!/_\chi G_v$ induced by the inclusion $R(\chi^i)\hookrightarrow\Gamma(\mathcal{O}_\Lambda)$ is a categorical quotient.
Moreover, there is an open subset $U$ of $\Lambda/\!/_\chi G_v$ such that $q^{-1}(U)=\Lambda_\chi^s$ and $q|_{\Lambda_\chi^s}:\Lambda_\chi^s\to U$ is a geometric quotient.
\end{prop}

Take the special character $\mathbf{1}=(1,\cdots,1)$. Then we can decide which points are $\mathbf{1}$-semistable or $\mathbf{1}$-stable by the following lemma (cf. [Nak2, 3.8] and [M2, Definition 4]). 

\begin{lem}\label{lem:sta}
One has $\Lambda_{\mathbf{1}}^{ss}=\Lambda_{\mathbf{1}}^s$.
A point $(A_i,B_i, \Delta_j,\Gamma_j)_{i,j}\in\Lambda$ is $\mathbf{1}$-semistable (or equivalently $\mathbf{1}$-stable) if and only if it satisfies the following condition :\\
For any linear subspaces $S_j$ of $V_j\,(j=1,\cdots,n)$ such that (1) $S_j$ contains $\mathrm{Im}\,\Gamma_j,\,(2) A_i(S_i)\\\subset S_{i+1}$ and (3) $B_i(S_{i+1})\subset S_i$ for all $i$ and $j$, one has $S_j=V_j$ for all $j$.
\end{lem}

We call a character $\chi$ of $G_v$ {\em generic} if $\Lambda_\chi^{ss}=\Lambda_\chi^s$. (For example, $\mathbf{1}\in\chi(G_v)$ is generic by the above lemma.) Quiver varieties have the following properties [Nak1, 4.1].

\begin{prop}
If a character $\chi$ is generic and $\mathfrak{M}_\chi(v,w)$ is nonempty, then $\mathfrak{M}_\chi(v,w)$ is a symplectic manifold and $\pi_\chi:\mathfrak{M}_\chi(v,w)\to\mathfrak{M}_\mathbf{0}(v,w)$ is a resolution of singularities of the image $\mathrm{Im}\,\pi_\chi$.
\end{prop}

\hspace{-6mm}\textbf{Remark.} \;The above proposition implies that $\mathrm{Im}\,\pi_\chi(\subset\mathfrak{M}_\mathbf{0}(v,w))$ is a symplectic variety and always has a crepant resolution. This holds for general quiver varieties.\vspace{3mm}

We can realize a nilpotent orbit closure $\overline{O_d}$ and its crepant resolution as quiver varieties. From now on we assume that $d\ne[1,\cdots,1]$ since $O_d$ is a point if $d=[1,\cdots,1]$. Let $a=[a_1\cdots a_m]$ be the dual partition of $d=[d_1,\cdots,d_l]$ i.e. the Young diagram corresponding to $a$ is obtained by transposing the Young diagram corresponding to $d$. For example, if $d=[4,4,2,1]$, then $a=[4,3,2,2]$. The corresponding Young diagrams are as follows.
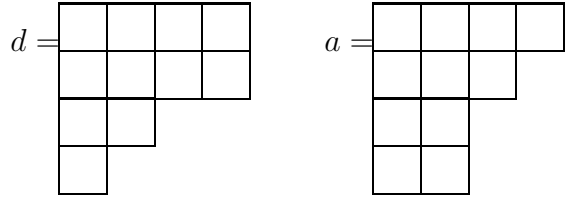
\vspace{3mm}

\begin{center}
$d=$\ydiagram{4,4,2,1}\hspace{10mm}$a=$\ydiagram{4,3,2,2}\\
\end{center}
Note that $m\ge2$ since $d\ne[1,\cdots,1]$.

Let $\mathcal{F}_a$ denote the flag variety of flag type $a$ i.e. 
$$\mathcal{F}_a=\{\mbox{subspaces }\{0\}=U_0\subset U_1\subset\cdots\subset U_{m-1}\subset U_m=\mathbf{C}^N|\mathrm{dim}\,U_i/U_{i-1}=a_i,\,1\le i\le m\}.$$
Then $G=SL_N(\mathbf{C})$ naturally acts on $\mathcal{F}_a$, and the stabilizer of the standard flag $F=(0\subset\mathbf{C}^{a_1}\subset\cdots\subset\mathbf{C}^{N-a_l}\subset\mathbf{C}^N)\in\mathcal{F}_a$ is a parabolic subgroup $P$ of $G$. Hence $\mathcal{F}_a$ is isomorphic to $G/P$. By [CG, 1.4], the cotangent bundle $T^*(G/P)$ is isomorphic to $G\times^P\mathfrak{n}(P)$. One can check that
$$G\times^P\mathfrak{n}(P)\to\{(x,U_\bullet)\in\overline{O_d}\times\mathcal{F}_a|x(U_{i+1})\subset U_i,\,0\le i\le m-1\}$$
$$\hspace{-57mm}(g,x)\mapsto(gxg^{-1},gF)$$
is an isomorphism. We identify $T^*\mathcal{F}_a$ with $\{(x,U_\bullet)\in\overline{O_d}\times\mathcal{F}_a|x(U_i)\subset U_{i+1},\,0\le i\le m-1\}$ via these isomorphisms. We define dimension vectors in $\mathbf{Z}^{m-1}$ as 
$$\tilde{v}=(N-a_1,N-a_1-a_2,\cdots,a_m)$$
and
$$\tilde{w}=(N,0,\cdots,0).$$
When we take $\mathbf{C}^N$ and $\mathbf{C}^{\tilde{v_i}}$'s as $W_1$ and $V_i$'s in the definition of $\Lambda$, an element of $\Lambda(\tilde{v},\tilde{w})$ corresponds to the following configuration.
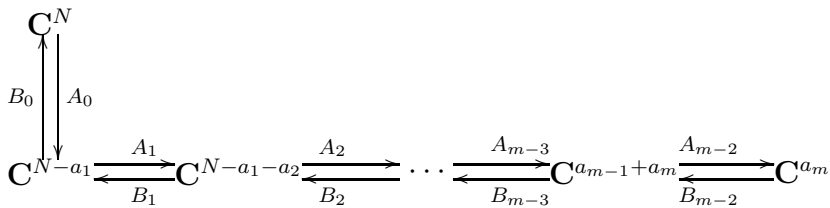
\vspace{5mm}

\begin{xy}
(0,0)*{}="",(15,0)*{\mathbf{C}^{N-a_1}}="1",(40,0)*{\mathbf{C}^{N-a_1-a_2}}="2",(65,0)*{\;\cdots\;}="dot",(90,0)*{\mathbf{C}^{a_{m-1}+a_m}}="n-1",(115,0)*{\mathbf{C}^{a_m}}="n",(15,20)*{\mathbf{C}^N}="w1"
\ar @<1mm>"1";"2" ^{A_1} \ar @<1mm> "2";"dot"^{A_2} \ar @<1mm> "dot";"n-1"^{A_{m-3}} \ar @<1mm> "n-1";"n"^{A_{m-2}}
\ar @<1mm> "2";"1"^{B_1} \ar @<1mm> "dot";"2"^{B_2} \ar @<1mm> "n-1";"dot"^{B_{m-3}} \ar @<1mm> "n";"n-1"^{B_{m-2}}
\ar @<1mm>"w1";"1"^{A_0} \ar @<1mm>"1";"w1"^{B_0}
\end{xy}
\vspace{5mm}

In the diagram we set $A_0:=\Gamma_1,B_0:=\Delta_1$ for convenience and omit $\{W_j\}_{2\le j\le m-1}$, $\{\Gamma_j\}_{1\le j\le m-1}$ and $\{\Delta_j\}_{1\le l\le m-1}$ because they are zero. For the special character $\tilde{\mathbf{1}}=(1,\cdots,1)\in\chi(G_{\tilde{v}})$, one easily sees that a point $(A_i,B_i)_{0\le i\le m-2}\in\Lambda(\tilde{v},\tilde{w})$ is $\tilde{\mathbf{1}}$-stable if and only if all $A_i$'s are surjective by Lemma \ref{lem:sta}.

\begin{prop}\label{prop:odbar}
[Nak1, \S7] One has the following commutative diagram.
$$\begin{CD} 
\mathfrak{M}_{\tilde{\mathbf{1}}}(\tilde{v},\tilde{w}) @>\sim>\theta>T^*\mathcal{F}_a \\ 
@V\pi_{\tilde{\mathbf{1}}}VV @VVpV \\ 
\mathfrak{M}_{\mathbf{0}}(\tilde{v},\tilde{w}) @>\sim>> \overline{O_d}      
\end{CD} $$
where $\theta$ is given by
$$(A_i,B_i)_i\mapsto (B_0A_0,(0\subset\mathrm{Ker}\,A_0\subset\mathrm{Ker}\,A_1A_0\subset\cdots\subset
\mathrm{Ker}\,A_{m-2}\cdots A_1A_0\subset\mathbf{C}^N))$$
and $p$ is the restriction of the first projection $\mathrm{End}(\mathbf{C}^N)\times\mathcal{F}_a\to \mathrm{End}(\mathbf{C}^N)$.
\end{prop}

\hspace{-6mm}{\bf Remark.} \;For a permutation $\sigma\in\mathfrak{S}_m$, we can define a new dimension vector $\sigma(\tilde{v})$ by replacing $a$ by $\sigma(a)=(a_{\sigma(1)},\cdots,a_{\sigma(m)})$ in the definition of $\tilde{v}$. One can show that the above proposition still holds if one replaces $\tilde{v}$, $a$ and $\mathfrak{M}_{\mathbf{0}}(\tilde{v},\tilde{w})$ by $\sigma(\tilde{v})$, $\sigma(a)$ and $\mathrm{Im}\,\pi_{\tilde{\mathbf{1}}}$ respectively. \vspace{3mm}

Next we consider the Slodowy slice $\mathcal{S}_{d',d}$ introduced in section \ref{sl} where $d'=[d'_1,\cdots,d'_n]$ be the partition of $N$ such that $O_{d'}\subset\overline{O_d}$. We define dimension vectors $v=(v_1,\cdots,v_{m-1})$ and $w=(w_1,\cdots,w_{m-1})$ in $\mathbf{Z}^{m-1}$ as
$$v_i=\sum_{k=i+1}^m(a'_k-a_k)\hspace{5mm} \mbox{and} \hspace{5mm}w_i=\sharp\{k|d'_k=i\}$$
where $a=[a_1\cdots a_m]$ and $a'=[a'_1\cdots a'_m]$ are the dual partitions of $d$ and $d'$ (now we allow $a'_k$'s to be zero).

Maffei proved that $\mathcal{S}_{d',d}$ and its crepant resolution are also obtained as quiver varieties.

\begin{thm}\label{thm:maf}([M2, Theorem 8] and its proof)
There are a injective homomorphism $\iota:G_v\hookrightarrow G_{\tilde{v}}$ and a closed immersion $\phi:\Lambda(v,w)\hookrightarrow\Lambda(\tilde{v},\tilde{w})$ such that:\\
(1) $\phi$ is $G_v$-equivariant.\\
(2) $\iota$ induces a surjective homomorphism $\psi:\chi(G_{\tilde{v}})\to\chi(G_v)$ such that $\psi(\tilde{\mathbf{1}})=\mathbf{1}$.\\
(3) $\phi$ induces closed immersions
$$\phi_1:\mathfrak{M}_{\mathbf{1}}(v,w)\hookrightarrow\mathfrak{M}_{\tilde{\mathbf{1}}}(\tilde{v},\tilde{w})\cong T^*\mathcal{F}_a$$
and
$$\phi_0:\mathfrak{M}_\mathbf{0}(v,w)\hookrightarrow\mathfrak{M}_\mathbf{0}(\tilde{v},\tilde{w})\cong\overline{O_d}$$
making the following diagram commutative:
$$\begin{CD} 
\mathfrak{M}_{\mathbf{1}}(v,w) @>\sim>\phi_1>\tilde{\mathcal{S}}_{d',d}\, @. :=p^{-1} (\mathcal{S}_x\cap\overline{O_d}) @. \subset T^*\mathcal{F}_a\\ 
@V\pi_{\mathbf{1}}VV @VV\mathrm{res}\,pV \\ 
\mathfrak{M}_{\mathbf{0}}(v,w) @>\sim>\phi_0> \mathcal{S}_{d',d} @. \cong \mathcal{S}_x\cap\overline{O_d}\subset\overline{O_d}
\end{CD} $$
where $x$ is some element of $O_{d'}$.
\end{thm}

\hspace{-6mm}\textbf{Remark.} 1.\;Just as Proposition \ref{prop:odbar}, the above claims still hold if one replaces $v$, $\tilde{v}$, $a$, $\mathfrak{M}_{\mathbf{0}}(v,w)$ and $\mathfrak{M}_{\mathbf{0}}(\tilde{v},\tilde{w})$ by $\sigma(v)$, $\sigma(\tilde{v})$, $\sigma(a)$, $\mathrm{Im}\,\pi_{\mathbf{1}}$ and $\mathrm{Im}\,\pi_{\tilde{\mathbf{1}}}$ respectively.

2.\;By (2), the homomorphism $\psi$ is an isomorphism if and only if $\chi(G_{\tilde{v}})_\mathbf{R}$ and $\chi(G_v)_\mathbf{R}$ have the same dimensions. The latter condition is also equivalent to the condition that every $v_i$ is nonzero.\\

Now we prepare some GIT terminologies for the next section. If two characters in $\chi(G_v)$ (resp. $\chi(G_{\tilde{v}})$) give the same semistable locus in $\Lambda(v,w)$ (resp. $\Lambda(\tilde{v},\tilde{w}))$ and hence give the  canonically isomorphic GIT quotients, we call them {\em GIT equivalent}. It is known that GIT equivalence classes give a chamber structure in $\chi(G_v)_\mathbf{R}$ (resp. $\chi(G_{\tilde{v}})_\mathbf{R}$) [T, 2.3] i.e.\\
(i) there are only finitely many GIT equivalence classes,\\
(ii) for every GIT equivalence class $C$, the closure $\overline{C}$ is a rational polyhedral cone in $\chi(G_v)_\mathbf{R}$ (resp. $\chi(G_{\tilde{v}})_\mathbf{R}$) and $C$ is a relative interior of $\overline{C}$.
\vspace{2mm}

\noindent We call $C$ a {\em GIT chamber} if $C$ is not contained in any hyperplane in $\chi(G_v)_\mathbf{R}$ (resp. $\chi(G_{\tilde{v}})_\mathbf{R}$). The complement of all GIT chambers in $\chi(G_v)_\mathbf{R}$ (resp. $\chi(G_{\tilde{v}})_\mathbf{R}$) is called a {\em wall}. It is known that $\chi\in\chi(G_{\tilde{v}})_\mathbf{R}$ (resp. $\chi(G_{\tilde{v}})_\mathbf{R}$) is generic if and only if $\chi$ is in a GIT chamber.

For the case of $\Lambda(\tilde{v},\tilde{w})$, we can describe the GIT chambers and the corresponding GIT quotients explicitly. The Weyl group $\mathfrak{S}_m$ naturally acts on $\chi(G_{\tilde{v}})_\mathbf{R}\cong\mathbf{R}^{m-1}$. This action coincides with the usual action of the Weyl group on the weight lattice if one changes the basis of $\mathbf{R}^{m-1}$ so that the standard basis $\mathbf{e}_1,\cdots,\mathbf{e}_{m-1}$ corresponds to the fundamental weights.
The GIT chambers coincide with the Weyl chambers $\{\sigma(C_\mathrm{fund})|\sigma\in\mathfrak{S}_m\}$ where $C_{\mathrm{fund}}=\mathbf{R}_{>0}^{m-1}\subset\mathbf{R}^{m-1}
\cong\chi(G_{\tilde{v}})_{\mathbf{R}}$
is the fundamental chamber which contains $\tilde{\mathbf{1}}$.

\begin{lem}
There is an algebraic isomorphism $\Phi_\sigma:\mathfrak{M}_\chi(\tilde{v},\tilde{w})\to\mathfrak{M}_{\sigma\chi}(\sigma(\tilde{v}),\tilde{w})$ over $\mathfrak{M}_{\mathbf{0}}(\tilde{v},\tilde{w})=\overline{O_d}$ for any $\sigma\in\mathfrak{S}_m$  and any generic $\chi\in\chi(G_{\tilde{v}})$ such that $\Phi_{\sigma\tau}=\Phi_\sigma\circ\Phi_\tau$ for all $\sigma,\tau\in\mathfrak{S}_m$.
\end{lem}
{\em Proof.} The isomorphism is constructed in [M1]. To check the compatibility with $\overline{O_d}$, it is sufficient to consider the case when $\sigma$ is the transposition $s_i=(i\;i+1)$ ($i=1,\cdots,m-1$) since every permutation $\sigma$ is generated by $s_i$'s. Note that $\tilde{v}_j=s_i(\tilde{v})_j$ for all $j\ne i$. For $x\in\mathfrak{M}_{\chi}(\tilde{v},\tilde{w})$, let $(A_i,B_i)_{0\le i\le m-2}\in\Lambda_{\chi}^{ss}(\tilde{v},\tilde{w})$ be any representative of $x$. By construction of $\Phi_{s_i}$, there is $(A'_i,B'_i)_{0\le i\le m-2}\in\Lambda_{s_i\chi}^{ss}(s_i(\tilde{v}),\tilde{w})$ such that\\
(1) $A_j=A'_j$ and $B_j=B'_j$ for all $j\ne i-1,i$.\\
(2) The sequence\\
$$0\to V'_i\xrightarrow{B'_{i-1}\oplus A'_i}V_{i-1}\oplus V_{i+1}\xrightarrow{A_{i-1}+B_i}V_i\to 0$$
is exact where $V_i$ (resp. $V'_i$) ($i=1,\cdots, m-1$) is the $\tilde{v}_i$ (resp. $s_i(\tilde{v})_i$)-dimensional vector space, $V_0=W_1=\mathbf{C}^N$ and $V_m=0$.\\
(3) The equality
$$(B'_{i-1}\oplus A'_i)\circ (A'_{i-1}+B'_i)=(B_{i-1}\oplus A_i)\circ (A_{i-1}+B_i)$$
in $\mathrm{End}\,(V_{i-1}\oplus V_{i+1})$ holds.\\
Then $\Phi_{s_i}(x)$ is defined to be the class of $(A'_i,B'_i)_{0\le i\le m-2}$ in $\mathfrak{M}_{s_i\chi}(s_i(\tilde{v}),\tilde{w})$. Since the image of $(A_i,B_i)_{0\le i\le m-2}$ (resp. $(A'_i,B'_i)_{0\le i\le m-2}$) in $\overline{O_d}$ is given by $B_0A_0$ (resp. $B'_0A'_0$) by Proposition \ref{prop:odbar} and its remark, the conditions (1) and (3) imply that $B_0A_0=B'_0A'_0$ and hence the compatibility with $\overline{O_d}$ follows.\qed

\vspace{3mm}
\hspace{-7mm}{\bf Remark.} \;The action of Weyl group on the quiver varieties is also defined in [Nak1], but this action is constructed analytically. \vspace{3mm}

By the above lemma, we see that $\mathfrak{M}_{\tilde{\chi}}(\tilde{v},\tilde{w})$ for a generic $\tilde{\chi}\in\chi(G_{\tilde{v}})$ is isomorphic to $T^*\mathcal{F}_{\sigma(a)}$ for some $\sigma\in\mathfrak{S}_m$. If $\sigma(a)\ne a$, the parabolic subgroup $P_{\sigma}$ of $G$ defined as the stabilizer of the standard flag in $\mathcal{F}_{\sigma(a)}$ is not conjugate to $P$. Moreover, every parabolic subgroup of $G$ whose Levi part is the same as $P$ is conjugate to $P_{\sigma}$ for some $\sigma$. Therefore the number of crepant resolutions of $\overline{O_d}$ is $\sharp\{\sigma(a)\in\mathbf{Z}_{>0}^m|\sigma\in\mathfrak{S}_m\}$ (cf. Theorem \ref{thm:nam}). The fact that every crepant resolution of $\overline{O_d}$ is of the form $T^*\mathcal{F}_{\sigma(a)}$ also follows from the argument in the next section.

To end this section, we give an example of a chamber structure on $\chi(G_{\tilde{v}})_{\mathbf{R}}$.

Suppose $d=[3,2,1]$. Then $a=[3,2,1]$ and hence $m=3$. In this case there are six GIT chambers on $\chi(G_{\tilde{v}})_{\mathbf{R}}\cong\mathbf{R}^2$ as in the figure below. To each GIT chamber, the corresponding GIT quotient is added there. 

\begin{center}
\setlength\unitlength{1truecm}
\begin{picture}(6,6)(0,0)
\put(3,3){\dashbox{0.2}(1,1)}
\put(0,3){\vector(1,0){6}}
\put(3,0){\vector(0,1){6}}
\put(5.5,0.5){\line(-1,1){5}}
\put(4,4){$\tilde{\mathbf{1}}$}
\put(4,2.5){$\mathbf{e}_1$}
\put(2.5,4){$\mathbf{e}_2$}
\put(4,5){$T^*\mathcal{F}_{(3,2,1)}$}
\put(1.3,5){$T^*\mathcal{F}_{(2,3,1)}$}
\put(0.3,3.7){$T^*\mathcal{F}_{(2,1,3)}$}
\put(1,1.5){$T^*\mathcal{F}_{(1,2,3)}$}
\put(3.2,1){$T^*\mathcal{F}_{(1,3,2)}$}
\put(4.5,2){$T^*\mathcal{F}_{(3,1,2)}$}
\put(2.5,2.5){$\textbf{O}$}
\end{picture}\\
\vspace{5mm}
Figure 1

\vspace{3mm}
The GIT chamber structure on $\chi(G_{\tilde{v}})_{\mathbf{R}}$ and the corresponding GIT quotients.
\end{center}

We will see later that $\chi(G_{\tilde{v}})_{\mathbf{R}}$ is naturally isomorphic to $H^2(T^*\mathcal{F}_{(3,2,1)},\mathbf{R})$ and the GIT chambers coincide with the ample cones via this isomorphism.

\section{The proofs of the theorems\label{proof}}
In this section we prove the theorems stated in section \ref{intro}. The key step is to identify some GIT chambers in the character groups with ample cones in Picard groups of crepant resolutions of $\overline{O_d}$ and $\mathcal{S}_{d',d}$.

Since $q:\Lambda_{\mathbf{1}}^{ss}(\tilde{v},\tilde{w})\to\mathfrak{M}_\mathbf{1}(\tilde{v},\tilde{w})$ is a geometric quotient by the connected group $G_{\tilde{v}}$, we have the following diagram of exact sequences of Abelian groups [KKV, 5.1].

$$\begin{CD} 
@.    1\\
@.    @VVV\\
@.    \mathcal{O}(\mathfrak{M}_{\tilde{\mathbf{1}}}(\tilde{v},\tilde{w}))^*/\mathbf{C}^*\\
@.    @VVV\\
@.    \mathcal{O}(\Lambda_{\tilde{\mathbf{1}}}^{ss}(\tilde{v},\tilde{w}))^*/\mathbf{C}^* @. 0 \\
@.     @VVV @VVV\\
@.     \chi(G_{\tilde{v}}) @. \mathrm{Pic}(\mathfrak{M}_{\tilde{\mathbf{1}}}(\tilde{v},\tilde{w}))\\
@.     @VVV @Vq^*VV\\
0@>>> \mathrm{Ker}\,h @>>> \mathrm{Pic}^{G_{\tilde{v}}}(\Lambda_{\tilde{\mathbf{1}}}^{ss}(\tilde{v},\tilde{w})) @>h>> \mathrm{Pic}(\Lambda_{\tilde{\mathbf{1}}}^{ss}(\tilde{v},\tilde{w}))\\ 
@.     @VVV @VVV \\
@.     0        @. \mathrm{torsion\;group}
\end{CD} $$
\vspace{3mm}
\\
where $\mathrm{Pic}^{G_{\tilde{v}}}$ means the Abelian group of $G_{\tilde{v}}$-equivariant isomorphism classes of line bundles with $G_{\tilde{v}}$-actions, and $h$ forgets $G_{\tilde{v}}$-actions. The image of $\tilde{\chi}\in\chi(G_{\tilde{v}})$ in $\mathrm{Pic}^{G_{\tilde{v}}}(\Lambda_{\tilde{\mathbf{1}}}^{ss}(\tilde{v},\tilde{w}))$ is the trivial line bundle with the twisted $G_{\tilde{v}}$-action by $\tilde{\chi}$. By the property of a quiver variety, the torsion group in the second column of the diagram is trivial [Nak1, 2.12]. So $q^*$ is an isomorphism and there is a map $f_{\tilde{\mathbf{1}}}:\chi(G_{\tilde{v}})\to \mathrm{Pic}(\mathfrak{M}_{\tilde{\mathbf{1}}}(\tilde{v},\tilde{w}))$ which commutes with the maps in the diagram. Note that in this case $R(\tilde{\chi})$ is equal to the vector space $H^0(\Lambda_{\tilde{\mathbf{1}}}^{ss}(\tilde{v},\tilde{w}), (q^*\circ f_{\tilde{\mathbf{1}}})(\tilde{\chi}))^{G_{\tilde{v}}}$ of $G_{\tilde{v}}$-invariant sections of $(q^*\circ f_{\tilde{\mathbf{1}}})(\tilde{\chi})$. In the same way we have $f_\mathbf{1}:\chi(G_v)\to\mathrm{Pic}(\mathfrak{M}_\mathbf{1}(v,w))$. This construction can be done for any generic $\chi\in\chi(G_v)$ and $\tilde{\chi}\in\chi(G_{\tilde{v}})$. So we also have $f_\chi:\chi(G_v)\to\mathrm{Pic}(\mathfrak{M}_\chi(v,w))$ and $f_{\tilde{\chi}}:\chi(G_{\tilde{v}})\to \mathrm{Pic}(\mathfrak{M}_{\tilde{\chi}}(\tilde{v},\tilde{w}))$.

The following proposition plays a key role in the proofs.
\begin{prop}\label{prop:comm}
The two maps $f_\mathbf{1}$ and $f_{\tilde{\mathbf{1}}}$ are isomorphisms and the following diagram commutes:
$$\begin{CD} 
\chi(G_{\tilde{v}})_\mathbf{R} @>\sim>f_{\tilde{\mathbf{1}}}>\mathrm{Pic}(\mathfrak{M}_{\tilde{\mathbf{1}}}(\tilde{v},\tilde{w}))_\mathbf{R} \\ 
@V\psi VV @VV\phi_\mathbf{1}^*V \\ 
\chi(G_v)_\mathbf{R} @>\sim>f_\mathbf{1}> \mathrm{Pic}(\mathfrak{M}_\mathbf{1}(v,w))_\mathbf{R}
\end{CD} $$
\end{prop}
{\em Proof.} The commutativity follows from the fact that $\phi:\Lambda_\mathbf{1}^{ss}(v,w)\hookrightarrow\Lambda_{\tilde{\mathbf{1}}}^{ss}(\tilde{v},\tilde{w})$ is $G_v$-equivariant. Note that $\mathfrak{M}_\mathbf{1}(v,w)$ (resp. $\mathfrak{M}_{\tilde{\mathbf{1}}}(\tilde{v},\tilde{w})$) is isomorphic to $\tilde{\mathcal{S}}_{d',d}$ (resp. $T^*\mathcal{F}_a$) by Theorem \ref{thm:maf}. As $\mathrm{dim}\,H^2(T^*\mathcal{F}_a,\mathbf{R})=\mathrm{dim}\,H^2(\mathcal{F}_a,\mathbf{R})=\mathrm{dim}\,\mathrm{Pic}(G/P)_{\mathbf{R}}$ is the dimension of the center of the Levi part of $P$ (cf. [Nam2, (P.3)]), it is equal to $m-1$, which is the dimesion of $\chi(G_{\tilde{v}})_\mathbf{R}$.  Consider the following commutative diagram:
$$\begin{CD} 
H^2(T^*\mathcal{F}_a,\mathbf{R}) @>>>H^2(\pi_\mathbf{1}^{-1}(x),\mathbf{R}) \\ 
@VVV @| \\ 
H^2(\tilde{\mathcal{S}}_{d',d},\mathbf{R}) @>>> H^2((\pi_\mathbf{1}|_{\tilde{\mathcal{S}}_{d',d}})^{-1}(x),\mathbf{R})
\end{CD} $$
where the horizontal arrows and the left vertical arrow are restriction maps. Since $\overline{O_d}$ (resp. $\mathcal{S}_{d',d}$) has a $\mathbf{C}^*$-acton with a unique fixed point [CG, 3.7], it is contractible. (This also follows from the fact that an affine quiver variety has a natural $\mathbf{C}^*$-acton with a unique fixed point; cf. [Nak1, \S5]) Hence the left vertical arrow is identified with $\phi_\mathbf{1}^*$ by Lemma \ref{lem:pic}. The upper horizontal map is surjective by [BO, Cororally 2.5], and the lower horizontal map is an isomorphism by [Nak1, 5.5]. Therefore $\phi_\mathbf{1}^*$ is surjective. The dimension $\mathrm{dim}\,\mathrm{Pic}(\mathfrak{M}_\mathbf{1}(v,w))_\mathbf{R}=\mathrm{dim}\,H^2(\pi_\mathbf{1}^{-1}(x),\mathbf{R})$ can be computed by the method in [BO, \S2] and one can check that it is equal to the number of nonzero $v_i$'s, which is also equal to $\mathrm{dim}\,\chi(G_v)_\mathbf{R}$. Thus the claim follows if we show that $f_{\tilde{\mathbf{1}}}$ is injective.

By considering the diagram at the beginning of this section, it is enough to show that $\mathcal{O}(\Lambda_{\tilde{\mathbf{1}}}^{ss}(\tilde{v},\tilde{w}))^*=\mathbf{C}^*$. Recall that
$$\Lambda_{\tilde{\mathbf{1}}}^{ss}(\tilde{v},\tilde{w})=\{(A_i,B_i)_{0\le i\le m-2}\in\Lambda_{\tilde{\mathbf{1}}}(\tilde{v},\tilde{w})|\mbox{all }A_i\mbox{'s are surjective}\}$$
by Lemma \ref{lem:sta}. Take any nowhere vanishing regular function $f\in\mathcal{O}(\Lambda_{\tilde{\mathbf{1}}}^{ss}(\tilde{v},\tilde{w}))^*$. We will show that $f$ is constant on $\Lambda_{\tilde{\mathbf{1}}}^{ss}(\tilde{v},\tilde{w})$. To this end, we think of $\Lambda_{\tilde{\mathbf{1}}}^{ss}(\tilde{v},\tilde{w})$ as a subset of the ``\{$A_i$\},\{$B_i$\}-plane'' and show that $f$ is constant both on the ``\{$A_i$\}-axis'' and along the ``\{$B_i$\}-directions''.

Since each $A_i$ is a nonsquare matrix, the subset $\Lambda_{\tilde{\mathbf{1}}}^{ss}(\tilde{v},\tilde{w})\cap\{B_0=\cdots=B_{m-2}=0\}$ is an affine space minus a codimension at least two locus. So the restriction of $f$ to this subset extends without poles to the affine space and hence $f$ is constant on this subset. For any surjective $A\in\bigoplus_{i=1}^{m-2}\mathrm{Hom}(V_i,V_{i+1})$, the subset $\Lambda_{\tilde{\mathbf{1}}}^{ss}(\tilde{v},\tilde{w})\cap\{(A_0,\cdots,A_{m-2})=A\}$ has a $\mathbf{C}^*$-action defined by
$$t:(A,B_i)_{0\le i\le m-2}\mapsto(A,tB_i)_{0\le i\le m-2}$$
for $t\in\mathbf{C}^*$. Note that this $\mathbf{C}^*$-action has a unique fixed point $(A,0)$. Therefore $\mathbf{C}^*$ acts on $f$ trivially (otherwise $f$ would have $(A,0)$ as a zero or a pole, which is not the case). This implies that $f$ is constant on $\Lambda_{\tilde{\mathbf{1}}}^{ss}(\tilde{v},\tilde{w})\cap\{(A_0,\cdots,A_n)=A\}$. We conclude that $f$ is constant on the whole $\Lambda_{\tilde{\mathbf{1}}}^{ss}(\tilde{v},\tilde{w})$. \qed

\vspace{5mm}
The following lemma shows the correspondence between the special GIT chamber in the character group and the ample cone in the Picard group. 

\begin{lem}\label{lem:ample}
Let $C_{\mathrm{fund}}$ (resp. $\tilde{C}_{\mathrm{fund}}$) be the GIT chamber in $\chi(G_v)_\mathbf{R}$ (resp. $\chi(G_{\tilde{v}})_\mathbf{R}$) which contains $\mathbf{1}$ (resp. $\tilde{\mathbf{1}}$). Then $f_\mathbf{1}(C_{\mathrm{fund}})$ (resp. $f_{\tilde{\mathbf{1}}}(\tilde{C}_{\mathrm{fund}}$)) is the $\pi_\mathbf{1}$-ample cone (resp. the $\pi_{\tilde{\mathbf{1}}}$-ample cone).
\end{lem}
{\em Proof.}
Let $\chi\in C_{\mathrm{fund}}$ and $L:=f_\mathbf{1}(\chi)$. Then
$$\begin{aligned}\mathrm{Proj}(\bigoplus_{i=0}^\infty H^0(\mathfrak{M}_\mathbf{1}(v,w),L^{\otimes i}))&=\mathrm{Proj}(\bigoplus_{i=0}^\infty H^0(\Lambda_\mathbf{1}^{ss}(v,w),q^*L^{\otimes i})^{G_v})\\&=\mathrm{Proj}(\bigoplus_{i=0}^\infty H^0(\Lambda_\chi^{ss}(v,w),q^*L^{\otimes i})^{G_v})\\&=\mathrm{Proj}\,(\bigoplus_{i=0}^\infty R(\chi^i))=\mathfrak{M}_\chi(v,w)=\mathfrak{M}_\mathbf{1}(v,w)\end{aligned}$$
and hence $L$ is $\pi_\mathbf{1}$-ample. Conversely, let $L\in\mathrm{Pic}(\mathfrak{M}_\mathbf{1}(v,w))_\mathbf{R}$ be $\pi_\mathbf{1}$-ample and $\chi:=f_{\mathbf{1}}^{-1}(L)\in C_{\mathrm{fund}}$. The line bundle $L$ gives a morphism
$$\phi_L:\mathfrak{M}_\mathbf{1}(v,w)\to\mathrm{Proj}(\bigoplus_{i=0}^\infty H^0(\mathfrak{M}_\mathbf{1}(v,w),L^{\otimes i}))$$
such that $\phi_L^*(\mathcal{O}(1))=L$. This is just the map $\Lambda/\!/_{\mathbf{1}} G_v\to\Lambda/\!/_\chi G_v$ of GIT quotients which is induced by the inclusion $\Lambda_\mathbf{1}^{ss}\cap\Lambda_\chi^{ss}\hookrightarrow\Lambda_\chi^{ss}$. Note that $\phi_L$ is an isomorphism since $L$ is $\pi_\mathbf{1}$-ample. This means that $\Lambda_\mathbf{1}^{ss}=\Lambda_\chi^{ss}$ and hence $\mathbf{1}$ and $\chi$ are GIT equivalent. For the case of $\tilde{C}_{\mathrm{fund}}$, the same proof works.
\qed
\vspace{5mm}

Take $L\in\mathrm{Pic}(\mathfrak{M}_\mathbf{1}(v,w))_\mathbf{R}$ from an ample cone in $\mathrm{Mov}(\pi_\mathbf{1})$ so generic that $\chi:=f_\mathbf{1}^{-1}(L)$ is not in a wall. Let $f:\mathfrak{M}_\mathbf{1}(v,w)\dashrightarrow Y$ be the corresponding rational map which is an isomorphism in codimension 1 over $\mathfrak{M}_{\mathbf{0}}(v,w)\cong\mathcal{S}_{d',d}$. We see that $Y$ can be obtained as the quiver variety associated to $\chi$ by the following lemma.

\begin{lem}\label{lem:mov}
For the $L$ and $\chi$ above, one has $\mathfrak{M}_{\chi}(v,w)=Y$ and $f_\mathbf{1}=f^*\circ f_\chi$ where $f^*:\mathrm{Pic}(Y)_\mathbf{R}\to\mathrm{Pic}(\mathfrak{M}_\mathbf{1}(v,w))_\mathbf{R}$ is the natural isomorphism (cf. Lemma \ref{lem:iso}).
\end{lem}
{\em Proof.} By the definition of $Y$ and the argument in the proof of the above lemma,
$$\begin{aligned}Y&=\mathrm{Proj}(\bigoplus_{i=0}^\infty H^0(\mathfrak{M}_\mathbf{1}(v,w),L^{\otimes i}))=\mathrm{Proj}(\bigoplus_{i=0}^\infty H^0(\Lambda_\mathbf{1}^{ss},q^*L^{\otimes i})^{G_v})\\
&=\mathrm{Proj}(\bigoplus_{i=0}^\infty H^0(\Lambda_\mathbf{1}^{ss}\cap\Lambda_\chi^{ss},q^*L^{\otimes i})^{G_v}).\end{aligned}$$
Note that $\mathfrak{M}_\mathbf{1}(v,w)$ and $\mathfrak{M}_{\chi}(v,w)$ are isomorphic in codimension 1 since they are both crepant resolutions of $\mathcal{S}_{d',d}$ (cf. Lemma \ref{lem:cre}). Since $\mathfrak{M}_\mathbf{1}(v,w)$ and $\mathfrak{M}_{\chi}(v,w)$ are the geometric quotients of $\Lambda_\mathbf{1}^{ss}$ and $\Lambda_\chi^{ss}$ respectively by Proposition \ref{prop:git}, $\Lambda_\mathbf{1}^{ss}$ and $\Lambda_\chi^{ss}$ are also isomorphic in codimension 1. Therefore
$$\begin{aligned}\mathfrak{M}_{\chi}(v,w)&=\mathrm{Proj}(\bigoplus_{i=0}^\infty H^0(\Lambda_\chi^{ss},q^*L^{\otimes i})^{G_v})\\&=\mathrm{Proj}(\bigoplus_{i=0}^\infty H^0(\Lambda_\mathbf{1}^{ss}\cap\Lambda_\chi^{ss},q^*L^{\otimes i})^{G_v})=Y.\end{aligned}$$
The latter claim follows directly from the construction of $f_\mathbf{1},f_\chi$ and $f^*$. \qed

\vspace{5mm}
\hspace{-6mm}\textbf{Remark.}\hspace{3mm}If we take a generic $\chi\in\chi(G_v)_\mathbf{R}$, then $Y:=\mathfrak{M}_{\chi}(v,w)$ is a crepant resolution of $\mathfrak{M}_{\mathbf{0}}(v,w)$ and there are two natural isomorphisms between relative Picard groups of $\mathfrak{M}_\mathbf{1}(v,w)$ and $Y$. The first is induced by the rational map $f:\mathfrak{M}_\mathbf{1}(v,w)\dashrightarrow Y$ which is isomorphic in codimension 1. The second is what comes from the fact that the relative Picard groups are isomorphic to $\chi(G_v)_\mathbf{R}$ via $f_\mathbf{1}$ and $f_\chi$. The above lemma implies that these two isomorphisms are the same if $L:=f_\mathbf{1}(\chi)$ is in the movable cone. But this fails if $L$ is not in $\mathrm{Mov}(\pi_\mathbf{1})$. Consider, for example, the cotangent bundle of the projective line. It is a unique crepant resolution of $\overline{O_{[2]}}$ and its Picard group is 1-dimensional. If $\chi=-\mathbf{1}$, then $L$ is anti-ample. The uniqueness of the crepant resolution implies that the rational map $f:\mathfrak{M}_\mathbf{1}(v,w)\dashrightarrow Y$ is an isomorphism over $\overline{O_{[2]}}$ and hence the first isomorphism $f^*$ preserves the ampleness. On the other hand, $f_{-\mathbf{1}}(-\mathbf{1})$ is ample (cf. Lemma \ref{lem:ample}) and therefore the second isomorphism $\mathrm{Pic}(Y)_\mathbf{R}\cong\mathrm{Pic}(\mathfrak{M}_\mathbf{1}(v,w))_\mathbf{R}$ sends an ample line bundle to an anti-ample one.

\vspace{5mm}
{\em Proof of Theorem \ref{thm:main1}}\\
Let $Y$ be a crepant resolution of $\overline{O_d}$. $Y$ is isomorphic to $\mathfrak{M}_{\tilde{\chi}}(\tilde{v},\tilde{w})$ over $\overline{O_d}$ for some generic $\tilde{\chi}\in\chi(G_{\tilde{v}})$ by Lemma \ref{lem:mov}. Then $\chi:=\psi(\tilde{\chi})$ is generic since $\psi$ is surjective.  The first claim of Theorem \ref{thm:main1} follows from the lemma below.

\begin{lem}
$\mathfrak{M}_\chi(v,w)$ is the restriction of $\mathfrak{M}_{\tilde{\chi}}(\tilde{v},\tilde{w})$ to $\mathcal{S}_{d',d}$.
\end{lem}
{\em Proof.} Let $Y'$ be the restriction of $Y$ to $\mathcal{S}_{d',d}$. We have
$$\phi^{-1}(\Lambda_{\tilde{\chi}}^{ss}(\tilde{v},\tilde{w}))\subset\Lambda_\chi^{ss}(v,w)$$
by the $G_v$-equivariance of $\phi$ (see Theorem \ref{thm:maf}). This inclusion induces an open immersion $Y'\to \mathfrak{M}_\chi(v,w)$ over $\mathcal{S}_{d',d}$ since $\mathfrak{M}_\chi(v,w)$ has the quotient topology determined by the geometric quotient $\Lambda_\chi^{ss}(v,w)\to\mathfrak{M}_\chi(v,w)$. It must be an isomorphism since $Y'\to\mathcal{S}_{d',d}$, which is the restriction of the projective morphism $\pi_{\tilde{\chi}}$ to $\mathcal{S}_{d',d}$, is proper. \qed 

\vspace{3mm}
Conversely let $Y'\to\mathcal{S}_{d',d}$ be any crepant resolution. Since $Y'$ and $\mathfrak{M}_\mathbf{1}(v,w)$ are isomorphic in codimension 1, there is an $L\in\mathrm{Pic}(\mathfrak{M}_\mathbf{1}(v,w))_\mathbf{R}$ such that the GIT quotient of $\Lambda(v,w)$ associated to $\chi:=f_\mathbf{1}^{-1}(L)$ is $Y'$ by Lemma \ref{lem:mov}. Take a generic $\tilde{\chi}$ such that $\psi(\tilde{\chi})=\chi$. Then $Y:=\mathfrak{M}_{\tilde{\chi}}(\tilde{v},\tilde{w})$ is a crepant resolution of $\overline{O_d}$. By the above lemma again, we see that $Y'=\mathfrak{M}_\chi(v,w)$ is the restriction of $Y$ to $\mathcal{S}_{d',d}$. This proves the second claim.

\vspace{5mm}
{\em Proof of Theorem \ref{thm:main2}}\\
The first claim follows from Lemma \ref{lem:pic}, Theorem \ref{thm:maf} and Proposition \ref{prop:comm}. The second claim follows from the lemma below.
 
\begin{lem}
If $\psi$ and $\phi_\mathbf{1}^*$ in Proposition \ref{prop:comm} are isomorphisms, then the set of ample cones in $\mathrm{Mov}(\pi_\mathbf{1})$ and the set of ample cones in $\mathrm{Mov}(\pi_{\tilde{\mathbf{1}}})$ are in one-to-one correspondence via $\phi_\mathbf{1}^*$.
\end{lem}
{\em Proof.} By Lemma \ref{lem:mov}, it is enough to show that $\phi_\mathbf{1}^*(\mathrm{Amp}(\pi_{\tilde{\mathbf{1}}}))=\mathrm{Amp}(\pi_\mathbf{1})$. The inclusion $\phi_\mathbf{1}^*(\mathrm{Amp}(\pi_{\tilde{\mathbf{1}}}))\subset\mathrm{Amp}(\pi_\mathbf{1})$ follows since $\phi_\mathbf{1}^*$ is the restriction map. To show the other inclusion, take $L$ from $\mathrm{Amp}(\pi_\mathbf{1})$. Set $\tilde{L}:=(\phi_\mathbf{1}^*)^{-1}(L)$. Note that we have an isomorphism
$$H_2(\mathfrak{M}_\mathbf{1}(v,w),\mathbf{R})\to H_2(\mathfrak{M}_{\tilde{\mathbf{1}}}(\tilde{v},\tilde{w}),\mathbf{R})\,$$
which is the dual of $\phi_\mathbf{1}^*$. Hence for every proper curve $\tilde{C}$ in $\mathfrak{M}_{\tilde{\mathbf{1}}}(\tilde{v},\tilde{w})$ which is contracted by $\pi_{\tilde{\mathbf{1}}}$, there is a proper curve $C$ in $\mathfrak{M}_\mathbf{1}(v,w)$ which is contracted by $\pi_\mathbf{1}$ and is homologically equivalent to $\tilde{C}$. Therefore $(\tilde{L}.\tilde{C})=(L.C)>0$ and $\tilde{L}$ is $\pi_{\tilde{\mathbf{1}}}$-ample. \qed

\vspace{5mm}
{\em Proof of Theorem \ref{thm:main3}}\\
To get the decomposition of $\mathcal{S}_{d',d}\,$, we should focus on $\Lambda(v,w)$ which corresponds to $\mathcal{S}_{d',d}$. If $v_i=0$ for some $i$, one easily sees by definition that $\Lambda(v,w)=\Lambda(v',w')\times\Lambda(v'',w'')$ where $v'=(v_1,\cdots,v_{i-1}),w'=(w_1,\cdots,w_{i-1}),v''=(v_{i+1},\cdots,v_{n})$ and $w''=(w_{i+1},\cdots,w_{n})$. By repeating this process, we get a decomposition
$$\Lambda(v,w)=\Lambda(v^1,w^1)\times\cdots\times\Lambda(v^r,w^r)$$ such that $v^i$ has no zero components for every $i$. Note that $G_v=\prod_{i=1}^r G_{v^i}$ and $\chi(G_v)=\prod_{i=1}^r\chi(G_{v^i})$. Since the action of $G_v$ on $\Lambda(v,w)$ equals the product of actions of $G_{v^i}$ on $\Lambda(v^i,w^i)$, the $\chi$-semistable locus $\Lambda_{\chi}^{ss}(v,w)$ equals the product $\prod_{i=1}^r\Lambda_{\chi_i}^{ss}$ for every $\chi\in\chi(G_v)$ where $\chi_i\in\chi(G_{v^i})$ is the $i$-th component of $\chi$. Therefore \begin{equation}
\Lambda(v,w)/\!/_\chi G_v=\prod_{i=1}^r\Lambda(v^i,w^i)/\!/_{\chi_i} G_{v^i}.
\end{equation}

Take $n_i$ so that $v^i$ is in $\mathbf{N}^{n_i}$. For each $i$, we define $\tilde{v}^i=(\tilde{v}_1^i,\cdots,\tilde{v}_{n_i}^i)$ and $\tilde{w}^i=(\tilde{w}_1^i,\cdots,\tilde{w}_{n_i}^i)$ as
$$\tilde{v}_j^i=\begin{cases}&v_j^i+\sum_{k=j+1}^{n_i}(k-j)w_k^i\hspace{20.1mm}\mbox{if }j=1,\cdots,n_i-1\\
&v_{n_i}^i\hspace{55mm}\mbox{if }j=n_i\end{cases}$$
and
$$ \tilde{w}_j^i=\begin{cases}&v_j^i+\sum_{k=1}^{n_i}kw_k^i\hspace{27.3mm}\mbox{if }j=1,\\
&0\hspace{52mm}\mbox{if }j=2,\cdots,n_i.\end{cases}$$
We also define partitions $d^i$ and $d'^i$ of $N_i:=\tilde{w}_1^i$ for each $i$ as
$$d_j^i=\left(\begin{aligned}&\mbox{the number of elements in }\{N^i-\tilde{v}_1^i,\tilde{v}_1^i-\tilde{v}_2^i,\cdots,\tilde{v}_{{n_i}-1}^i-\tilde{v}_{n_i}^i\}\\
&\mbox{ which are at least } j\end{aligned}\right)$$
and 
$$d'^i=\left(\begin{aligned}&\mbox{the Young diagram such that the number of rows which} \\&\mbox{consist of }j\mbox{ boxes is }w_j^i.\end{aligned}\right)$$
By Theorem \ref{thm:maf} and its remark, one can check that there is a closed immersion $\mathrm{Im}\,\pi_{\mathbf{1}}\hookrightarrow\mathrm{Im}\,\pi_{\tilde{\mathbf{1}}}$ which is isomorphic to $\mathcal{S}_{d'^i,d^i}\hookrightarrow\overline{O_{d^i}}$ where $\mathrm{Im}\,\pi_{\mathbf{1}}$ (resp. $\mathrm{Im}\,\pi_{\tilde{\mathbf{1}}}$) is the image of $\pi_{\mathbf{1}}:\mathfrak{M}_\mathbf{1}(v^i,w^i)\to\mathfrak{M}_\mathbf{0}(v^i,w^i)$ (resp. $\pi_{\tilde{\mathbf{1}}}:\mathfrak{M}_{\tilde{\mathbf{1}}}(\tilde{v}^i,\tilde{w}^i)\to\mathfrak{M}_\mathbf{0}(\tilde{v}^i,\tilde{w}^i)$).
Claim (i) in Theorem \ref{thm:main3} follows from the equation (5.1), since every crepant resolution is obtained as $\mathfrak{M}_{\chi}(v,w)$ for some $\chi$ by Lemma \ref{lem:mov}. The restriction map $H^2(X_i,\mathbf{R})\to H^2(Y_i,\mathbf{R})$ in claim (ii) is an isomorphism from Proposition \ref{prop:comm} since $\chi(G_{v^i})_\mathbf{R}$ and $\chi(G_{\tilde{v}^i})_\mathbf{R}$ have the same dimensions.

For claim (iii), the decomposition $\chi(G_v)=\prod_{i=1}^r\chi(G_{v^i})$ gives a decomposition of the cohomology group by Proposition \ref{prop:comm}. Let $Z'_i$ be another crepant resolution of $\mathcal{S}_{d'^i,d^i}$ for $i=1,\cdots,r$ and $Z'=Z_1'\times\cdots\times Z_r'$. Note that $Y'$ and $Z'$ are isomorphic in codimension one over $\mathcal{S}_{d',d}$ if and only if $Y_i'$ and $Z'_i$ are isomorphic in codimension one over $\mathcal{S}_{d'^i,d^i}$ for all $i$. Since GIT chambers in $\chi(G_v)_\mathbf{R}$ are the products of chambers in $\chi(G_{v^i})_\mathbf{R}$'s, the ample cones in $H^2(Y',\mathbf{R})$ are the products of ample cones in $H^2(Y_i',\mathbf{R})$'s. This completes the proof of Theorem \ref{thm:main3}.

\section{Examples\label{eg}}
In this section we apply the main results to concrete examples. Let $d$ and $d'$ be partitions of a natural number $N$ such that $O_{d'}\subset\overline{O_d}$. One can find $d'^i$'s and $d^i$'s in the decomposition of $\mathcal{S}_{d',d}$ in Theorem \ref{thm:main3} by following its proof. By using Young diagrams, this can be done visually. Let $a=[a_1,\cdots,a_m]$ (resp. $a'=[a'_1,\cdots,a'_{m'}]$) be the dual partition of $d$ (resp. $d'$). For natural numbers $p\le q(\le m)$, we define a new Young diagram $d_{p,q}$ (resp. $d'_{p,q}$) whose $i$-th column consists of $a_{i+p-1}$ (resp. $a'_{i+p-1}$) boxes for $1\le i\le q-p+1$.
\vspace{3mm}

For example, if $d=$\ydiagram{5,4,2}\;, then $d_{2,4}$=\ydiagram{3,3,1}.
\vspace{3mm}

The Young diagrams corresponding to $d'^i$'s and $d^i$'s are obtained as follows:
\vspace{2mm}

{\bf Step 1.} Divide $d$ (resp. $d'$) into diagrams of the forms of $d_{p,q}$ (resp. $d'_{p,q}$) such that (the number of boxes in $d_{p,q}$)=(the number of boxes in $d'_{p,q}$) as fine as possible.
\vspace{2mm}

{\bf Step 2.} Remove the first common rows in $d_{p,q}$ and $d'_{p,q}$ for each $p,q$.
\vspace{2mm}

Then the resulting partitions are the desired ones.
\vspace{3mm}

{\bf Example 1.} \;The case where $d=[5,4,3,3,2,1]$ and $d'=[4,4,4,2,2,1,1]$.
\vspace{1mm}

\begin{center}
$d=$\ydiagram{5,4,3,3,2,1}\hspace{10mm}$d'=$\ydiagram{4,4,4,2,2,1,1}
\end{center}
\vspace{7mm}

{\bf Step 1} \hspace{5mm} These Young diagrams are divided as follows.
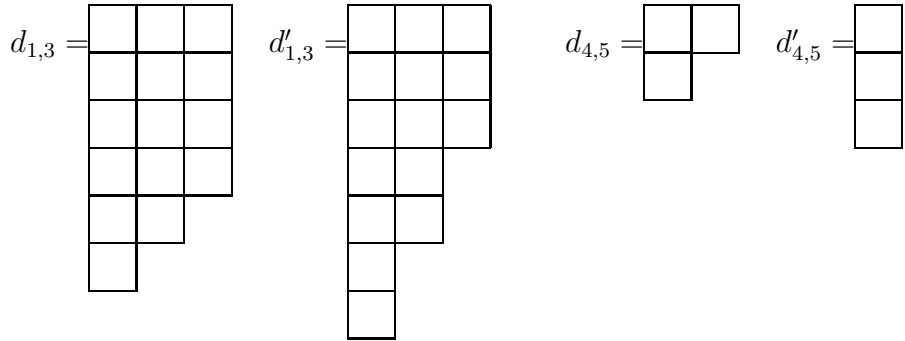
\vspace{3mm}

\begin{center}
$d_{1,3}=$\ydiagram{3,3,3,3,2,1}\hspace{5mm}$d'_{1,3}=$\ydiagram{3,3,3,2,2,1,1}
\hspace{10mm}$d_{4,5}=$\ydiagram{2,1}\hspace{5mm}$d'_{4,5}=$\ydiagram{1,1,1}\\
\end{center}
\vspace{7mm}

{\bf Step 2} \hspace{5mm} By removing the first common rows, we get
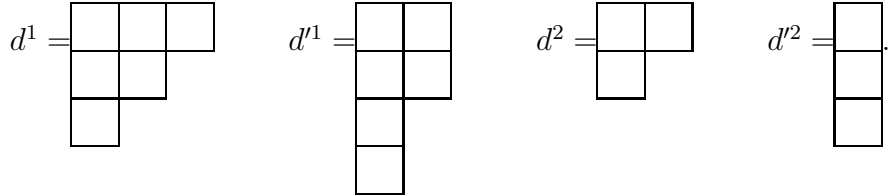
\begin{center}
$d^1=$\ydiagram{3,2,1}\hspace{10mm}$d'^1=$\ydiagram{2,2,1,1}\hspace{10mm}
$d^2=$\ydiagram{2,1}\hspace{10mm}$d'^2=$\ydiagram{1,1,1}.\\\vspace{3mm}
\end{center}
\vspace{3mm}

Hence $\mathcal{S}_{d',d}\cong\mathcal{S}_{d'^1,d^1}\times\mathcal{S}_{d'^2,d^2}$.
Since the numbers of crepant resolutions of $\overline{O_{d^1}}$ and $\overline{O_{d^2}}$ are
$\sharp\{\sigma(3,2,1)|\sigma\in\mathfrak{S}_3\}=3!$ and $\sharp\{\sigma(2,1)|\sigma\in\mathfrak{S}_2\}=2!$ respectively, the number of crepant resolutions of $\mathcal{S}_{d',d}$ is $3!\times2!=12$.
\vspace{7mm}

{\bf  Example 2.} \;The case where $d=[5,4,3,1]$ and $d'=[5,3,3,2]$.
\vspace{3mm}

\begin{center}
$d=$\ydiagram{5,4,3,1}\hspace{10mm}$d'=$\ydiagram{5,3,3,2}
\end{center}
\vspace{8mm}

{\bf Step 1} \hspace{5mm} These Young diagrams are divided as follows.
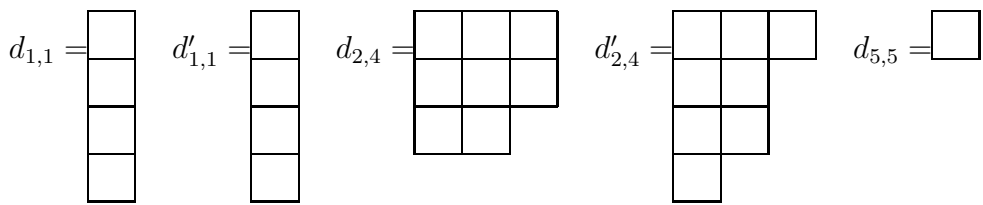
\vspace{3mm}

\begin{center}
$d_{1,1}=$\ydiagram{1,1,1,1}\hspace{5mm}$d'_{1,1}=$\ydiagram{1,1,1,1}
\hspace{5mm}$d_{2,4}=$\ydiagram{3,3,2}\hspace{5mm}$d'_{2,4}=$\ydiagram{3,2,2,1}
\hspace{5mm}$d_{5,5}=$\ydiagram{1}
\end{center}
\vspace{2mm}

$d'_{5,5}=$\ydiagram{1}
\vspace{5mm}

{\bf Step 2} \hspace{5mm} By removing the first common rows, we get
\vspace{3mm}

\begin{center}
$d^1=$\ydiagram{3,2}\hspace{10mm}$d'^1=$\ydiagram{2,2,1}.
\end{center}
\vspace{3mm}

Hence $\mathcal{S}_{d',d}\cong\mathcal{S}_{d'^1,d^1}$.
Since the numbers of crepant resolutions of $\overline{O_{d^1}}$ is $\sharp\{\sigma(2,2,1)|\sigma\in\mathfrak{S}_3\}=3$, that of $\mathcal{S}_{d',d}$ is also 3. Note that the number of crepant resolutions of $\overline{O_d}$ is $\sharp\{\sigma(4,3,3,2,1)|\sigma\in\mathfrak{S}_5\}=60$, which is much larger than that of $\mathcal{S}_{d',d}$.


\begin{thebibliography}{99}
\bibitem[B]{Beauville}
A. Beauville,
Symplectic singularities, Invent. Math. {\bf 139} (2000), 541-549.

\bibitem[BCHM]{BirkarCasciniHaconMcKernan}
C. Birkar, P. Cascini, C., Hacon and J. McKernan, 
Existence of minimal models for varieties of log general type, J. Amer. Math. Soc. {\bf 23} (2010), 405-468.

\bibitem[BO]{BrundanOstrik}
J. Brundan, V. Ostrik,
Cohomology of Spaltenstein varieties, Transform. Groups {\bf 16} (2011), 619-648.

\bibitem[CG]{ChrissGinzburg}
N. Chriss and V. Ginzburg,
Representation theory and complex geometry, Birkh\"{a}user, Boston, 1997.

\bibitem[CM]{CollingwoodMcGovern}
D. Collingwood, and W. McGovern,
Nilpotent orbits in semi-simple Lie algebras, van Nostrand Reinhold, Math. Series, 1993.

\bibitem[DH]{DolgachevHu}
I. Dolgachev and Y. Hu,
Variation of geometric invariant theory quotients, Inst. Hautes Etudes Sci. Publ. Math. (1998), no. {\bf 87}, 5-56.

\bibitem[HK]{HuKeel}
Y. Hu and S. Keel,
Mori dream spaces and GIT, Michigan Math. J. {\bf 48} (2000),  331-348.

\bibitem[KKV]{KnopKraftVust}
F. Knop, H. Kraft and T. Vust,
The Picard group of a G-variety, Algebraische
Transformationsgruppen und Invariantentheorie, DMV Semin. {\bf 13}, 77-87 (1989).

\bibitem[KP]{KraftProcesi}
H. Kraft and C. Procesi,
Closures of conjugacy classes of matrices are normal, Invent. Math. {\bf 53} (1979), 227-247.

\bibitem[LNS]{LehnNamikawaSorger}
M. Lehn, Y. Namikawa and C. Sorger,
Slodowy slices and universal Poisson deformations, Compos. Math. {\bf 148} (2012), 121-144.

\bibitem[M1]{Maffei1}
A. Maffei,
A remark on quiver varieties and Weyl groups. Ann. Scuola Norm. Sup. Pisa Cl. Sci. (5), {\bf 1}(2002), 649-686.

\bibitem[M2]{Maffei2}
A. Maffei,
Quiver varieties of type A, Comment. Math. Helv. {\bf 80} (2005), 1-27.

\bibitem[MFK]{MumfordFogartyKirwan}
D. Mumford, J. Fogarty and F. Kirwan,
Geometric invariant theory. Springer-Verlag, Berlin, 1994.

\bibitem[Nak1]{Nakajima1}
H. Nakajima,
Instantons on ALE spaces, quiver varieties, and Kac-Moody algebras. Duke Math. J. {\bf 76} (1994), 365-416.

\bibitem[Nak2]{Nakajima2}
H. Nakajima,
Quiver varieties and Kac-Moody algebras. Duke Math. J. {\bf 91} (1998), 515-560.

\bibitem[Nam1]{Namikawa1}
Y. Namikawa,
Birational geometry of symplectic resolutions of nilpotent orbits. Moduli spaces and arithmetic geometry, 75–116, Adv. Stud. Pure Math., {\bf 45}, Math. Soc. Japan, Tokyo, 2006. 

\bibitem[Nam2]{Namikawa2}
Y. Namikawa,
Birational geometry and deformations of nilpotent orbits, Duke Math. J. {\bf 143} (2008), 375-405.

\bibitem[Nam3]{Namikawa3}
Y. Namikawa, Poisson deformations and birational geometry, J. Math. Sci. Univ. Tokyo {\bf 22} (2015), no.~1, 339--359. 

\bibitem[P]{Panyushev}
D.I. Panyushev,
Rationality of singularities and the Gorenstein properties of nilpotent orbits. Functional Anal. Appl. {\bf 25}, 225-226(1991).

\bibitem[S]{Slodowy}
P. Slodowy,
Simple singularities and simple algebraic groups, Lecture Notes in Mathem., {\bf 815}, SpringerVerlag, 1980.

\bibitem[T]{Thaddeus}
M. Thaddeus,
Geometric invariant theory and flips, J. Amer. Math. Soc. {\bf 9} (1996), 691-723.

\end{thebibliography}
\end{document}